\documentclass[AMS,STIX1COL]{WileyNJD-v2}

\usepackage{moreverb}
\usepackage{Local_Definitions}
\usepackage[capitalise]{cleveref}

\usepackage{siunitx}
\usepackage{pgfplots}
\usepackage{pgfplotstable}
\usepackage{booktabs}
\usepackage{xcolor}
\hypersetup{
  colorlinks   = true, %Colours links instead of ugly boxes
  urlcolor     = blue, %Colour for external hyperlinks
  linkcolor    = blue, %Colour of internal links
  citecolor   =  purple%Colour of citations
}

\usepackage{paralist}

\usepackage{todonotes}
\usepackage{subfiles}
\usepackage{tikz}
\usetikzlibrary{fit,positioning,arrows,calc,backgrounds}
\newtheorem{example}[theorem]{Example}

\newcommand\BibTeX{{\rmfamily B\kern-.05em \textsc{i\kern-.025em b}\kern-.08em
T\kern-.1667em\lower.7ex\hbox{E}\kern-.125emX}}

\articletype{Original Paper}%

\received{<day> <Month>, <year>}
\revised{<day> <Month>, <year>}
\accepted{<day> <Month>, <year>}

\definecolor{darkgreen}{rgb}{0,0.6,0}
\definecolor{darkblue}{rgb}{0,0,0.6}
\definecolor{darkred}{rgb}{0.8,0.0,0.0}
\usepackage{subcaption}
\tikzset{%
    L1abserror/.style={darkgreen, solid, mark size=0pt},
    L2abserror/.style={darkblue, solid, mark size=0pt},
    Linfabserror/.style={darkred, solid, mark size=0pt},
    L1relerror/.style={darkgreen, solid, mark size=0pt},
    L2relerror/.style={darkblue, solid, mark size=0pt},
    Linfrelerror/.style={darkred, solid, mark size=0pt}
    }

\begin{document}

\title{Three Ways to Solve Partial Differential Equations with Neural Networks --- A Review%
\protect\thanks{This work was supported by the German Federal Ministry Education and Research (BMBF) as part of the project SOPRANN -- Synthese optimaler Regelungen und adaptiver Neuronaler Netze für Mobilitätsanwendungen (05M20OCA)}}

\author[1]{Jan Blechschmidt}
\author[1]{Oliver G. Ernst}
\authormark{Blechschmidt, Ernst}
\address[1]{\orgdiv{Department of Mathematics}, \orgname{TU Chemnitz}, \orgaddress{\state{Saxony}, \country{Germany}}}
\corres{Jan Blechschmidt, Reichenhainer Str.\ 41, 09126 Chemnitz \email{jan.blechschmidt@math.tu-chemnitz.de}}
\presentaddress{Jan Blechschmidt, Reichenhainer Str.\ 41, 09126 Chemnitz}

\abstract[Abstract]{%
Neural networks are increasingly used to construct numerical solution methods for partial differential equations.
In this expository review, we introduce and contrast three important recent approaches attractive in their simplicity and their suitability for high-dimensional problems: physics-informed neural networks, methods based on the Feynman-Kac formula and methods based on the solution of backward stochastic differential equations.
The article is accompanied by a suite of expository software in the form of Jupyter notebooks in which each basic methodology is explained step by step, allowing for a quick assimilation and experimentation.
An extensive bibliography summarizes the state of the art.
}


\jnlcitation{\cname{\author{J.\ Blechschmidt}, \author{O.~G.\ Ernst}}
    (\cyear{2021}),
\ctitle{Three ways to solve linear and nonlinear partial differential equations using neural networks --- A Review}, \cjournal{GAMM Mitteilungen 2021}, \cvol{TODO}.}

\maketitle
\section{Introduction} \label{sec:Intro}

	The spectacular successes of neural networks in machine learning tasks such as computer vision, natural speech processing and game theory as well as the prospect of harnessing the computing power of specialized hardware such as Google's Tensor Processing Units and Apple's Neural Engine designed to efficiently execute neural networks has led the scientific community to investigate their suitability also for high performance computing tasks.
	The result is now an exciting new research field known as \emph{scientific machine learning}, where techniques such as deep neural networks and statistical learning are applied to classical problems of applied mathematics.
		In this expository survey our intention is to provide an accessible introduction to recent developments in the field of numerical solution of linear and nonlinear partial differential equations (PDEs) using techniques from machine learning and artificial intelligence.
		
	After decades of research on the numerical solution of PDEs, manifold challenges remain.
	One that applies to essentially all classical discretization schemes is that they suffer from the \emph{curse of dimensionality} first formulated by Bellman in the 1950s in the context of optimal control problems \cite{Bellman1957}.
	In its simplest manifestation (see \cite{Powell2009} for a more extensive discussion) this notion states that doubling of the number of degrees of freedom in each of $d$ coordinate directions increases the solution complexity (at least) by a factor of $2^d$.
	In a similar spirit, the number of degrees of freedom when discretizing a 100-dimensional PDE with only 10 nodes in each coordinate direction exceeds the estimated number of atoms in the universe (around $10^{80}$) by several orders of magnitude.
	One might think that equations in such high dimensions have little practical relevance, but they are common in mathematical finance and portfolio optimization where the spatial dimension is determined by the number of financial assets in the market.
	Other areas prone to high-dimensional PDE problems include stochastic control, differential games and quantum physics.
	The challenge of solving high-dimensional PDEs has been taken up in a number of papers, and are addressed in particular in \cref{sec:FeynmanKacSolver} for linear Kolmogorov PDEs and in \cref{sec:DeepBSDESolver} for semilinear PDEs in non-divergence form.
	Another impetus for the development of data-driven solution methods is the effort often necessary to develop tailored solution methods for different kinds of nonlinear PDEs.
This will play a particular role in \cref{sec:PINNs}.

Neural networks offer attractive approximation capabilities for highly nonlinear functions. 
Their compositional nature contrasts with the more conventional additive form of trial functions in linear function spaces in which PDE solution approximations are constructed by Galerkin, collocation or finite volume methods.
Their computational parametrization through statistical learning and large-scale optimization methods using modern hardware, software systems and algorithms are making them increasingly amenable for solving nonlinear and high-dimensional PDEs.

PDE solvers based on (deep) neural networks typically cannot compete with classical numerical solution methods in low to moderate dimensions -- in particular as solving an algebraic equation is generally simpler than solving the highly nonlinear large-scale optimization problems associated with neural network training.
	Moreover, they currently lack the mature error analysis that has been established for traditional numerical methods.
	In addition, many specialized methods have been developed over the years for specific problems, often incorporating constraints or physical assumptions directly into the approximations.
	On the other hand, the ease with which methods such as the \emph{physics-informed neural networks} to be discussed below can be applied to essentially any differential equation makes them attractive for rapid prototyping when efficiency and high accuracy are not the principal concern.

% While we aim to give a useful overview, the publication this area is incredibly high.
While we aim to provide a useful overview, research activity in this area is incredibly intense and impossible to cover exhaustively.
Therefore, we have decided to present three approaches that have generated a lot of interest in recent years in detail in Sections \ref{sec:PINNs}--\ref{sec:DeepBSDESolver}.
Further scientific machine learning methods for solving PDEs are collected in~\cref{sec:FurtherMethods}.
Additionally, we want to draw some attention to another recent overview~\cite{BeckOverview} which contains many references, in particular works focusing on the solution of PDEs in high-dimensions.

A unique feature of this paper is a collection of accompanying Jupyter notebooks that contain sample Python implementations of the methods reviewed in~\cref{sec:PINNs,sec:FeynmanKacSolver,sec:DeepBSDESolver} with detailed comments and explanations as well as a number of numerical experiments.
The notebooks are freely available from the GitHub repository \url{https://github.com/janblechschmidt/PDEsByNNs} and can even be executed in \emph{Google Colaboratory} directly in a web browser with no need for local installations.
Of course, the reader may also download and run the notebooks on her local machine.

The remainder of the paper is organized as follows: \cref{sec:PINNs} discusses physics-informed neural networks, a straightforward and flexible approach for leveraging machine learning technology on challenging nonlinear PDE problems. 
\cref{sec:FeynmanKacSolver,sec:DeepBSDESolver} are devoted to recent methods based on the long-established link between PDEs and stochastic processes, which for high dimensions makes approximations based on sampling attractive due to their  dimension independence. 
Here neural networks on dedicated hardware can make the sample-based training very efficient.
\cref{sec:FurtherMethods} provides an outlook to related developments in this area followed by a concluding \cref{sec:conclusion}.

%-----------------------------------------------------------------------------------------
\section{Physics-Informed Neural Networks} \label{sec:PINNs}

The flexibility of deep neural networks as a universal technique for function approximation comes at the price of a large number of parameters to be determined in the supervised learning phase, and therefore typically demands a large volume of training data.
Physics-informed neural networks (PINNs) are a scientific machine learning technique for solving partial differential equation (PDE) problems in the \emph{small data} setting, meaning only the PDE problem data is available rather than a large number of value pairs of the indepenent and dependent variables.
PINNs generate approximate solutions to PDEs by training a neural network to minimize a loss function consisting of terms representing the misfit of the initial and boundary conditions along the boundary of the space-time domain as well as the PDE residual at selected points in the interior.
While precursors of this approach date back to the early 1990s
\cite{%
LeeKang1990,%
PsichogiosUngar1992,%
LagarisEtAl1998,%
LagarisEtAl2000},
the term PINN as well as a surge of ensuing research activity was initiated by the two-part report  
\cite{raissi2017physicsI,raissi2017physicsII} subsequently published in \cite{RaissiPerdikarisKarniadakis2019}.

We describe the PINN approach for approximating the solution $u:[0,T] \times \domain \to \R$ of an evolution equation
\begin{subequations} \label{eq:PINN_IBVP}
\begin{align} \label{eq:PINN_NonlinearPDE}
    \partial_t u (t,x) + \NN[u](t,x) &= 0, && (t,x) \in (0,T] \times \domain, \\
    u(0,x) &= u_0(x), \quad && x \in \domain,
\end{align}
where $\NN$ is a nonlinear differential operator acting on $u$, 
$\domain \subset \R^d$ a bounded domain,
$T$ denotes the final time and
$u_0: \domain \to \R$ the prescribed initial data.
Although the methodology allows for different types of boundary conditions, we restrict our discussion to the inhomogeneous Dirichlet case and prescribe
\begin{align} \label{eq:PINNBC}
        \hspace{7em} u(t,x) &= u_b(t,x),  && \quad (t,x) \in (0,T] \times \partial \domain,
\end{align}
\end{subequations}
where $\partial \domain$ denotes the boundary of the domain $\domain$ and $u_b: (0,T] \times \partial \domain \to \R$ the given boundary data.
The method constructs a neural network approximation $u_\theta(t,x) \approx u(t,x) $ of the solution of \eqref{eq:PINN_IBVP}, where $u_\theta :[0,T] \times \domain \to \R$ denotes a function realized by a neural network with parameters $\theta$.

In contrast to other learning-based methods that try to infer the solution by a purely data-driven approach, i.e., by fitting a neural network to a number of state-value pairs $\{(t_i,x_i,u(t_i,x_i))\}_{i=1}^N$, PINNs take the underlying PDE (the ``physics'') into account.
Taking advantage of modern machine learning software environments, which provide automatic differentiation capabilities for functions realized by neural networks, the approximate solution $u_\theta$ is differentiated with respect to the time and space variables, which allows the residual of the nonlinear PDE~\eqref{eq:PINN_NonlinearPDE} to be evaluated at a set of collocation points. 
In this way, the physics encoded in the differential equation is made available for a loss function measuring the extent to which the PDE problem \eqref{eq:PINN_IBVP} is satisfied by $u_\theta$.

While the focus of other methods employing neural networks for solving PDEs is on mitigating the curse of dimensionality in high dimensions, the strength of PINNs lies in their flexibility in that they can be applied to a great variety of challenging PDEs, whereas classical numerical approximations typically require tailoring to the specifics of a particular PDE.
In particular, this includes problems from computational physics that are notoriously hard to solve with classical numerical approaches due to, e.g., strong nonlinearities, convection dominance or shocks, see also the last paragraph in~\cref{sec:PINNExtension}.
A further challenge that can be addressed by this approach is the regime with a small number of data samples, which is common for physical experiments since the acquisition of new data samples is often expensive.

In \cite{raissi2017physicsI} the authors introduce the PINN methodology for solving nonlinear PDEs and demonstrate its efficiency for the Schrödinger, Burgers and Allen-Cahn equations.
The focus of the second part \cite{raissi2017physicsII} lies in the simultaneous solution of a nonlinear PDE of the form~\eqref{eq:PINN_NonlinearPDE} and the identification of corresponding unknown parameters $\lambda$ which enter the nonlinear part of the differential equation.
This problem setting has been studied within the regime of Gaussian processes  in \cite{RaissiKarniadakis2018,RaissiPerdikarisKarniadakis2017GP,Rudye1602614}.
For both problem settings, the authors discuss, depending on the type of data available, a time-continuous and time-discrete approach.
We discuss these methods next.

\subsection{Continuous Time Approach} \label{sec:PINNCont}

The continuous time approach for the parabolic PDE~\eqref{eq:PINN_IBVP} as described in~\cite{raissi2017physicsI} is based on the (strong) residual of a given neural network approximation $u_\theta \colon [0,T] \times \domain \to \R $ of the solution $u$ with respect to \eqref{eq:PINN_NonlinearPDE}
\begin{align} \label{eq:PINN}
    r_\theta (t,x) := \partial_t u_\theta (t,x) + \NN[u_\theta] (t,x).
\end{align}
The neural network class considered here are \emph{multilayer feed-forward neural networks}, sometimes known as  \emph{multilayer perceptrons}.
Such networks are compositions of alternating affine linear $W^\ell \cdot + b^\ell $ and nonlinear functions $\sigma^\ell(\cdot)$ called activations, i.e.,
\begin{align*}
    u_\theta(z) 
    := 
    W^L \sigma^L(W^{L-1} \sigma^{L-1} ( \cdots \sigma^1 (W^0 z + b^0) \cdots ) + b^{L-1} ) + b^L,
\end{align*}
where $W^\ell$ and $b^\ell$ are weight matrices and bias vectors, and $z=[t,x]^T$.
This highly nonlinear compositional structure of the approximating function $u_\theta$ forms the core of many neural network-based machine learning methods, and has been found to possess remarkably good approximation properties in many applications.

In general, training a neural network, i.e., determining the (typically large number of) parameters $\theta$, using gradient-based optimization methods 
\cite{%
glorot2010understanding,%
goodfellow2016deep,%
ruder2016overview,%
BottouEtAl2018} 
such as stochastic gradient descent \cite{BottouEtAl2018}, the Adam optimizer \cite{kingma2014adam}, or AdaGrad \cite{duchi2011adaptive}, requires the derivative of $u_\theta$ with respect to its unknown parameters $W^\ell$ and $b^\ell$.
To incorporate the PDE residual \eqref{eq:PINN} into the loss function to be minimized, PINNs require a further differentiation to evaluate the differential operators $\partial_t u_\theta$ and $\NN[u_\theta]$.
Thus the PINN term $r_\theta$ shares the same parameters as the original network $u_\theta(t,x)$, but respects the ``physics'' of~\eqref{eq:PINN_NonlinearPDE}.
Both types of derivatives can be easily obtained by automatic differentiation \cite{baydin2018automatic} with current state-of-the-art machine learning libraries, e.g., TensorFlow \cite{tensorflow2015-whitepaper} or PyTorch \cite{paszke2017automatic}.
In \cref{rem:PINN_Example} below, we show how such a PINN can be derived explicitly for the one-dimensional time-dependent eikonal equation.

The PINN approach for the solution of the PDE~\eqref{eq:PINN_IBVP} now proceeds by minimization of the loss functional
\begin{align} \label{eq:PINNLossFunctional}
    \phi_\theta(X) := \phi_\theta^r(X^r) + \phi_\theta^0(X^0) + \phi_\theta^b(X^b),
\end{align}
where $X$ denotes the collection of training data and the loss function $\phi_\theta$ contains the following terms:
\begin{itemize}
\item the mean squared residual
\begin{align*}
   \phi_\theta^r(X^r) := \frac{1}{N_r}\sum_{i=1}^{N_r} \left|r_\theta\left(t_i^r, x_i^r\right)\right|^2
\end{align*}
in a number of collocation points $X^r:=\{(t_i^r, x_i^r)\}_{i=1}^{N_r} \subset (0,T] \times \domain$, where $r_\theta$ is the physics-informed neural network~\eqref{eq:PINN},
\item the mean squared misfit with respect to the initial and boundary conditions
\begin{align*}
   \phi_\theta^0(X^0) 
   := 
   \frac{1}{N_0}
   \sum_{i=1}^{N_0} \left|u_\theta\left(t_i^0, x_i^0\right) - u_0\left(x_i^0\right)\right|^2
   \quad \text{ and } \quad
   \phi_\theta^b(X^b) 
   := 
   \frac{1}{N_b}
   \sum_{i=1}^{N_b} \left|u_\theta\left(t_i^b, x_i^b\right) - u_b\left(t_i^b, x_i^b\right)\right|^2
\end{align*}
in a number of points $X^0:=\{(t^0_i,x^0_i)\}_{i=1}^{N_0} \subset \{0\} \times \domain$ and $X^b:=\{(t^b_i,x^b_i)\}_{i=1}^{N_b} \subset (0,T] \times \partial \domain$, where $u_\theta$ is the neural network approximation of the solution $u\colon[0,T] \times \domain \to \R$.
\end{itemize}
We note that the training data $X$ consists entirely of time-space coordinates. 
Moreover, individual weighting  of each loss term in~\eqref{eq:PINNLossFunctional} may help improve the convergence of the scheme, see e.g.\ \cite{rao2020flow}.

\subsubsection{Example: Burgers Equation} \label{sec:ExampleBurgers}

To illustrate the PINN approach we consider the one-dimensional Burgers equation on the spatial domain $\DD = [-1,1]$
\begin{equation} \label{eq:BurgersEQ}
\begin{aligned}
    \partial_t u + u \, \partial_x u - (0.01/\pi) \, \partial_{xx} u &= 0, \quad &&\quad (t,x) \in (0,1] \times (-1,1),\\
   u(0,x) &= - \sin(\pi \, x),                 \quad &&\quad x \in [-1,1],\\
   u(t,-1) = u(t,1) &= 0,                      \quad &&\quad t \in (0,1].
\end{aligned}
\end{equation}
This PDE arises in various disciplines such as traffic flow, fluid mechanics and gas dynamics, and can be derived from the Navier-Stokes equations, see \cite{BasdevantDevilleHaldenwangLacroixOuazzaniPeyretOrlandi1986}.
We assume that the collocation points $X^r$ as well as the points for the initial and boundary data $X^0$ and $X^b$ are generated by random sampling from a uniform distribution.
Although uniformly distributed data are sufficient in our experiments, the authors of \cite{raissi2017physicsI} employed a space-filling Latin hypercube sampling strategy \cite{Stein1987}.
Our numerical experiments indicate that this strategy slightly improves the observed convergence rate, but for simplicity the code examples accompanying this paper employ uniform sampling throughout.

We choose training data of size $N_0 = N_b =\num{50}$ and $N_r = \num{10000}$.
In this example, adopted from \cite{raissi2017physicsI}, we assume a deep neural network of the following structure:
the input is scaled elementwise to lie in the interval $[-1,1]$, followed by 8 fully connected layers each containing 20 neurons and each followed by a hyperbolic tangent activation function and one output layer.
This setting results in a network containing $\num{3021}$ trainable parameters (first hidden layer: $2 \cdot 20 + 20 = 60$; seven intermediate layers: each $20 \cdot 20 + 20 = 420$; output layer: $20 \cdot 1 + 1 = 21$).

The loss functional~\eqref{eq:PINNLossFunctional} can be minimized by a number of algorithms, our accompanying code implements gradient descent-based algorithms as well as a variant of the limited-memory Broyden--Fletcher--Goldfarb--Shanno (BFGS) algorithm \cite{liu1989limited} which was also used in the numerical experiments in~\cite{raissi2017physicsI}.
Although currently the majority of neural networks are trained with gradient descent-based methods, BFGS is a quasi-Newton algorithm also often employed for scientific machine learning tasks.

\begin{figure}[htpb]
\centering
\begin{minipage}{.48\textwidth}
   \includegraphics[width=0.98\linewidth]{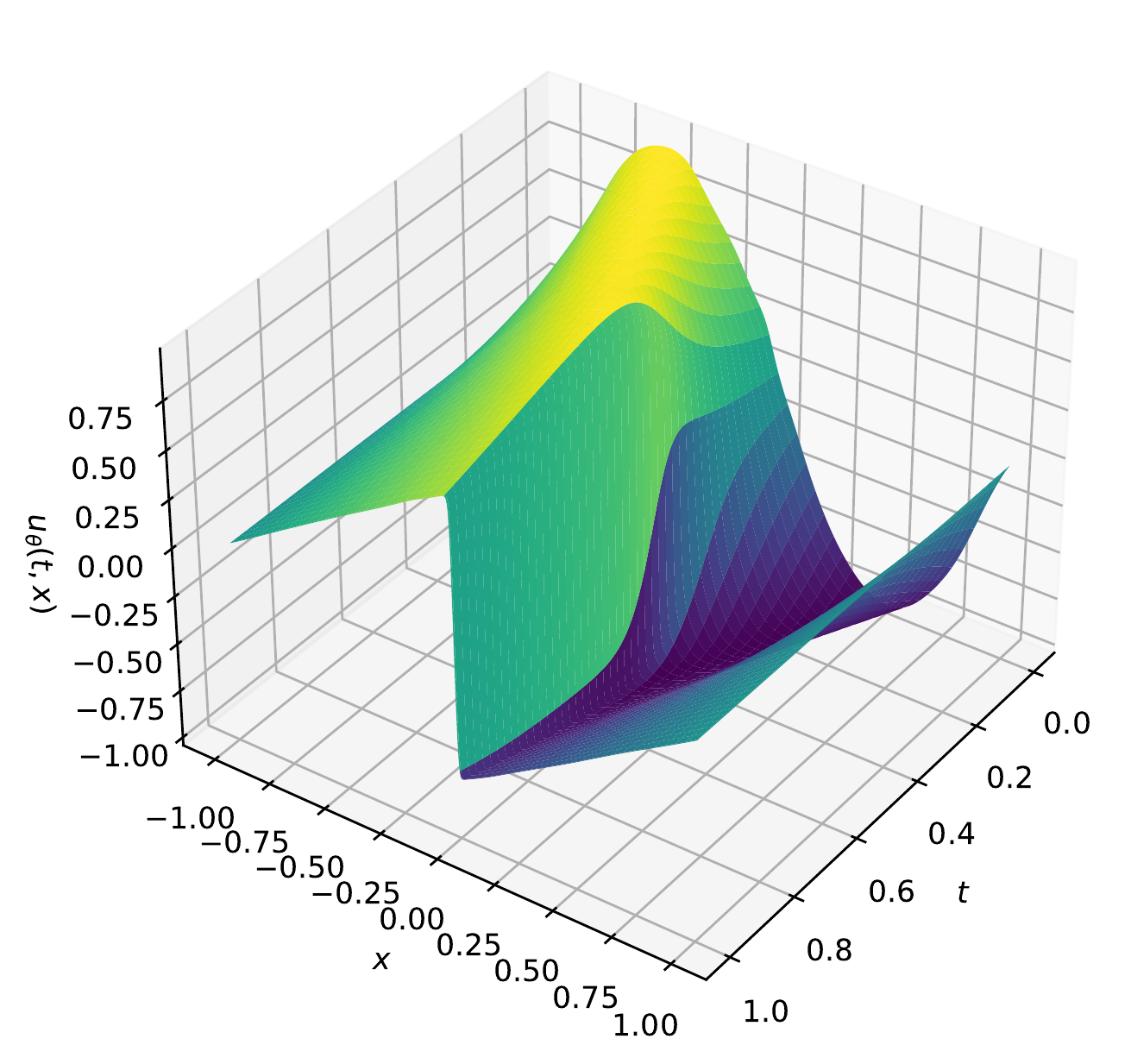}
\end{minipage}%
\begin{minipage}{.48\textwidth}
   \includegraphics[width=0.98\linewidth]{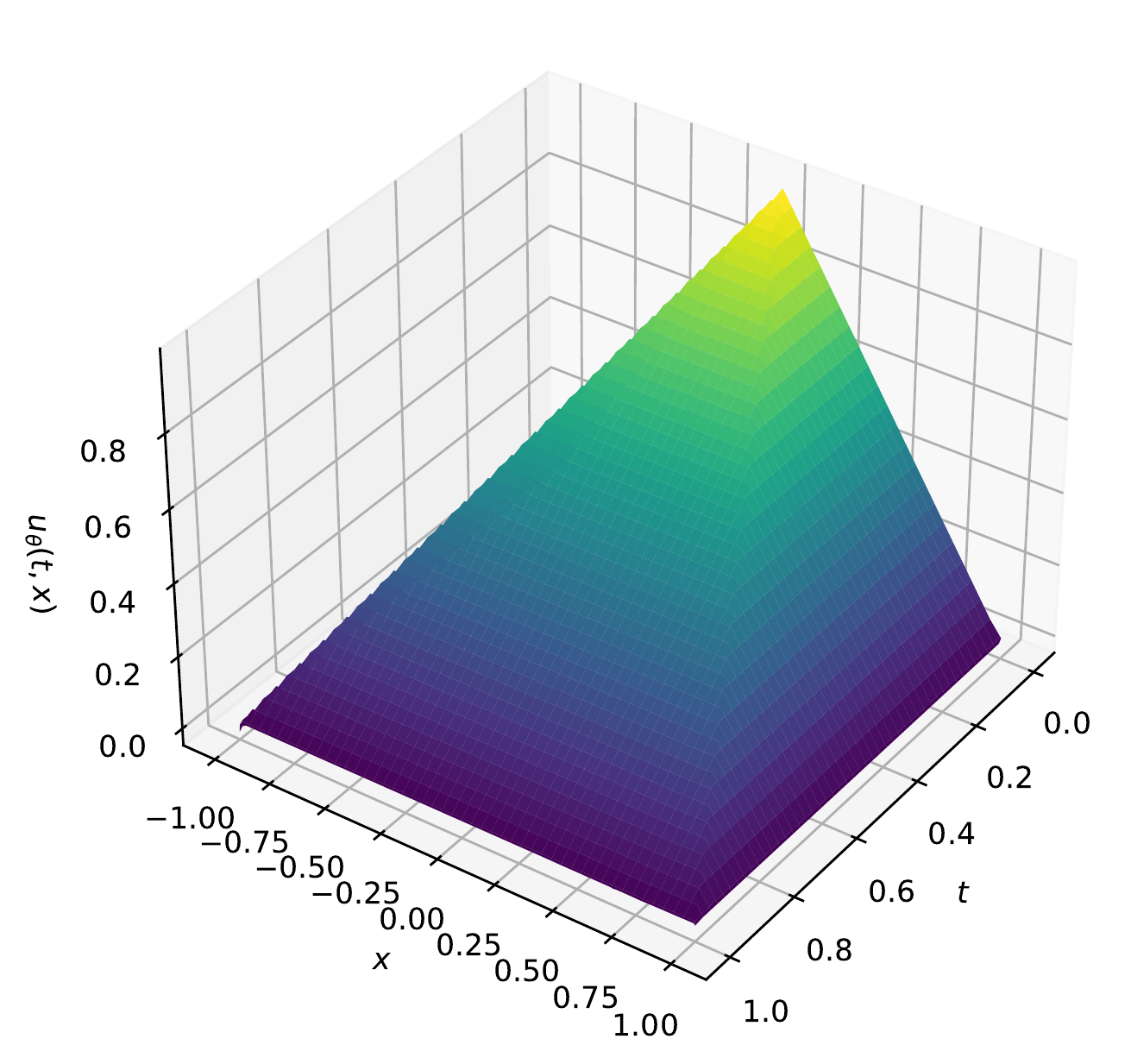}
\end{minipage}
\caption{\textbf{Left}: PINN approximation $u_\theta$ of the solution of Burgers equation~\eqref{eq:BurgersEQ}. 
The shock formation at around $t=0.4$ is clearly visible. 
\textbf{Right}: Approximate solution $u_\theta$ of the eikonal equation~\eqref{eq:Eikonal} with sharp edges at $t = \abs{x}$. Both examples are implememented in the accompanying Jupyter Notebook \texttt{PINN\_Solver.ipynb}.}%
    \label{fig:BurgersEQ}
\end{figure}

The left panel of \cref{fig:BurgersEQ} shows the approximate solution of the Burgers equation~\eqref{eq:BurgersEQ} after $\num{5000}$ training epochs with the Adam optimizer and learning rate\footnote{The chosen learning rates used in the Adam optimizer in this section are not based on any hyperparameter optimization but were selected in a way that ensured stable and reliable results.}
$\delta(n) = 0.01 \, \textbf{1}_{\{n < \num{1000}\}} + 0.001 \, \textbf{1}_{\{\num{1000} \le n < \num{3000}\}} + \num{0.0005} \, \textbf{1}_{\{\num{3000} \le n\}}$ which decays in a piecewise constant fashion.

\subsubsection{Example: Eikonal Equation} \label{sec:EikonalEquation}

As a second example we consider the one-dimensional eikonal equation backward in time on the domain $\DD=[-1,1]$
\begin{align}
\begin{aligned} \label{eq:Eikonal}
   -\partial_t u(t,x) + \abs{\nabla u}(t,x) &= 1, 
   		\quad & &(t,x) \in [0,T) \times [-1,1],\\
   u(T,x) &= 0, \quad & &x \in [-1,1],\\
   u(t,-1) = u(t, 1) &= 0, \quad & & t \in [0,T).
\end{aligned}
\end{align}
Note that the partial differential equation in~\eqref{eq:Eikonal} can be equally written as a Hamilton-Jacobi-Bellman equation, viz
\begin{align*}
    -\partial_t u(t,x) + \sup_{\abs{c} \le 1} \{c \, \nabla u(t,x)\} = 1, \qquad (t,x) \in [0,T) \times [-1,1],
\end{align*}
which characterizes the solution of an optimal control problem seeking to minimize the distance from a point $(t,x)$ to the boundary $[0,T] \times \partial \domain \cup \{T\} \times \domain$.
As is easily verified, the solution is given by $u(t,x) = \min\{ 1 - t, 1 - \abs{x} \}$.
The fact that~\eqref{eq:Eikonal} runs backward in time is in accordance with its interpretation as the optimality condition of a control problem.
Note that~\eqref{eq:Eikonal} is transformed into a forward evolution problem~\eqref{eq:PINN_NonlinearPDE} by the change of variables $\hat t = T - t$.

The neural network model chosen for this particular problem can be simpler.
We decided to use only two hidden layers with 20 neurons in each, resulting in $\num{501}$ unknown parameters (first hidden layer: $2 \cdot 20 + 20 = 60$; one intermediate layer: $20 \cdot 20 + 20 = 420$; output layer: $20 \cdot 1 + 1 = 21$).
To account for the lack of smoothness of the solution, we choose a non-differentiable activation function, although the hyperbolic tangent function seems to be able to approximate the kinks in the solution sufficiently well.
Here, we decided to use the \emph{leaky rectified linear unit (leaky ReLU)} activation function \cite{MaasEtAl2013}
\begin{align*}
    \sigma(z) = \begin{cases}
        z &\text{ if } z \ge 0,\\
        0.1 \, z &\text{ otherwise},
    \end{cases}
\end{align*}
which displays a non-vanishing gradient when the unit is not active, i.e., when $z < 0$.
The approximate solution after $N_\text{epochs} =\num{10000}$ epochs of training with the Adam optimizer \cite{kingma2014adam} and a piecewise constant learning rate 
\begin{equation} \label{eq:lrEikonal}
    \delta(n) = 0.1 \textbf{1}_{\{n < \num{3000}\}} +  0.01 \, \textbf{1}_{\{\num{3000} \le n < \num{7000}\}} +  0.001 \, \textbf{1}_{\{\num{7000} \le n\}}
\end{equation}
is displayed in the right panel of \cref{fig:BurgersEQ}.
    Noting that the explicit solution of the eikonal equation is a piecewise linear function on a convex polyhedral domain, closer inspection yields the closed-form expression
\[
	u(t,x) = \relu(x+1) - \relu(x+t) - \relu(x-t),
\]
which can be represented exactly by a neural network with one hidden layer containing three neurons.
    In order to study the capability of the PINN approach combined with the Adam optimizer to recover the solution of this problem we conducted an experiment for which we counted the number of successful attempts to train the model to achieve a training loss below the threshold $\phi_\theta(X) < 10^{-10}$.
    Otherwise, when a maximum number of iterations of $\num{100000}$ was reached, the algorithm had most often converged to a local minimum and no further decrease of the loss could be expected.
    We compared the activation functions leaky $\relu$ (slope $0.1$ for negative values) and standard $\relu$ (zero slope for negative values) on a set of different network architectures for ten uniformly drawn sets of training data with $N_r = \num{2000}$, $N_0 = 25$ and $N_b=50$ with learning rate as given in~\eqref{eq:lrEikonal}.
\Cref{table:Eikonal} shows the absolute number of successes among ten independent runs, indicating clearly that the leaky $\relu$ outperforms standard $\relu$ in this case.

\pgfplotstableset{
    multicolumn names, % allows to have multicolumn names
    col sep=comma, % the seperator in our .csv file
    sci zerofill,
    columns/col/.style={
        column name=Activation,
        string type},
    columns/l1_n3/.style={
        column name=Dim,
        int detect},
    columns/l1_n10/.style={
        column name=Dim,
        int detect},
    columns/l1_n25/.style={
        column name=Dim,
        int detect},
    columns/l2_n3/.style={
        column name=Dim,
        int detect},
    columns/l2_n10/.style={
        column name=Dim,
        int detect},
    columns/l2_n25/.style={
        column name=Dim,
        int detect},
    every head row/.style={
        output empty row},
    every last row/.style={after row=\bottomrule}, % rule at bottom
    every head row/.append style={before row={
        \toprule
        & \multicolumn{3}{c}{One hidden layer} &  \multicolumn{3}{c}{Two hidden layers}\\
        Activation & 3 Neurons & 10 Neurons & 25 Neurons & 3 Neurons & 10 Neurons & 25 Neurons\\
    \midrule}},
    % every head row/.style={
    %     before row={\toprule}, % have a rule at top
    %     after row={\midrule}, % rule under units
    %     % Activation & $\nneurons=3$ & $\nneurons=10$ & $\nneurons=25$ & $\nneurons=3$ & $\nneurons=10$ & $\nneurons=25$\\}
    %     % \\
    % },
}
    % every head row/.append style={before row={
    %     & \multicolumn{3}{c}{Layers $1$} & \multicolumn{3}{c}{Layers $2$}\\
    %     Activation & d1 & d2 & d3 & d11 & d12 & d13}
    % }},
\begin{table}[htb]
    \begin{center}
        \begin{filecontents*}{./table_Eikonal.csv}
col,l1_n3,l1_n10,l1_n25,l2_n3,l2_n10,l2_n25
ReLU,0, 3, 8, 0, 5, 3
Leaky ReLU,0, 9, 10, 1, 7, 10
\end{filecontents*}
\pgfplotstabletypeset[
            columns={{col},{l1_n3},{l1_n10},{l1_n25},{l2_n3},{l2_n10},{l2_n25}},
        ]{./table_Eikonal.csv}
    \end{center}
    \caption{Number of successful attempts to learn the solution of the eikonal equation \eqref{eq:Eikonal} for different network architectures for ten randomly initialized sets of training data with $N_r = \num{2000}$, $N_0 = 25$ and $N_b=50$.
    An attempt is considered successful if it achieves a training loss below the threshold $\phi_\theta(X) < 10^{-10}$.}
    \label{table:Eikonal}
\end{table}

\begin{figure}[htpb]
\centering
\tikzset{%
    neuron/.style={
        circle,
        draw,
        fill=yellow,
        fill opacity=.2,
        text opacity=1,
        minimum size=1cm
    },
    pinn_neuron/.style={
        circle,
        draw,
        fill=blue,
        fill opacity=.2,
        text opacity=1,
        minimum size=1cm
    },
    headline/.style={
        minimum size=1cm
    },
    subhead/.style={
        minimum size=1cm,
        font=\footnotesize
    },
    background/.style={rectangle,
    fill=gray!20,
    draw,
    inner sep=0.5cm,
    rounded corners=5mm},
    varbackground/.style={rectangle,
    fill=gray!10,
    inner sep=0.2cm,
    rounded corners=5mm}
    }
    \begin{tikzpicture}[x=1.5cm, y=1.5cm, >=stealth]
        \node [headline] (input-head) at (0, 3.8) {Input layer};
        \node [headline] (h-1-head) at (2, 3.8) {Hidden layer};
        \node [headline] (output-head) at (4, 3.8) {Output layer};
        \node [headline] (derivative) at (7.0, 3.8) {Derivative layer};
        \node [headline] (output-pinn) at (9, 3.8) {PINN};

        \node [subhead] (h-1-subhead) at (2, 3.50) {$z = W [t,x]^T + b$};
        \node [subhead] (output-1-subhead) at (4, 3.50) {$u_\theta = U\,\sigma(z) + c$};
        \node [subhead] (derivative-subhead) at (7, 3.50) {$\nabla u_\theta = (U \, \diag(\sigma'(z)) \, W)^T$};
        \node [subhead] (output-pinn-subhead) at (9, 3.50) {$r_\theta = - \tfrac{\partial u_\theta}{\partial t} + \left|\tfrac{\partial u_\theta}{\partial x}\right| - 1$};

        \node [neuron] (input-1) at (0,2.0) {$t$};
        \node [neuron] (input-2) at (0,1.0) {$x$};
        \node [neuron] (h-1-1) at (2,2.5) {$z_1$};
        \node [neuron] (h-1-2) at (2,1.5) {$z_2$};
        \node [neuron] (h-1-3) at (2,0.5) {$z_3$};
        \node [neuron] (output-1) at (4,1.5) {$u_\theta$};
        \node [pinn_neuron] (output-2) at (7,2.0) {$\frac{\partial u_\theta}{\partial t}$};
        \node [pinn_neuron] (output-3) at (7,1.0) {$\frac{\partial u_\theta}{\partial x}$};
        \node [pinn_neuron] (pinn) at (9,1.5) {$r_\theta$};

        \foreach \i in {1,2}
        \foreach \j in {1,2,3}
        \draw [->] (input-\i) -- (h-1-\j);

        \foreach \i in {1,2,3}
        \foreach \j in {1}
        \draw [->] (h-1-\i) -- (output-\j);

        % \foreach \i in {2,3}
        % \draw [->] (output-1) -- (output-\i);
        \foreach \i in {2,3}
        \draw [->] ($(4,1.5)+(1.0cm,0)$) -- (output-\i);

        \foreach \i in {2,3}
        \draw [->] (output-\i) -- (pinn);

    \begin{pgfonlayer}{background}
        \node [background,
                    fit=(input-1) (h-1-1) (h-1-3) (output-1)] {};
        % \node [varbackground,
        %             fit=(output-2) (output-3) (pinn)] {};
        % \node [varbackground,
        %             fit=(u_0) (grad_0)] {};
        % \node [background,
        %             fit=(grad_1) (u_Nm1)] {};
        % \node [background,
        %             fit=(u_N)] {};
    \end{pgfonlayer}

    \end{tikzpicture}
\caption{Illustration of a neural network with a single hidden layer (yellow). 
Complete network includes the physics-informed neural network $r_\theta$ for the one-dimensional eikonal equation~\eqref{eq:Eikonal} derived from the spatial and temporal derivatives of $u_\theta$.}%
    \label{fig:SimpleNN}
\end{figure}
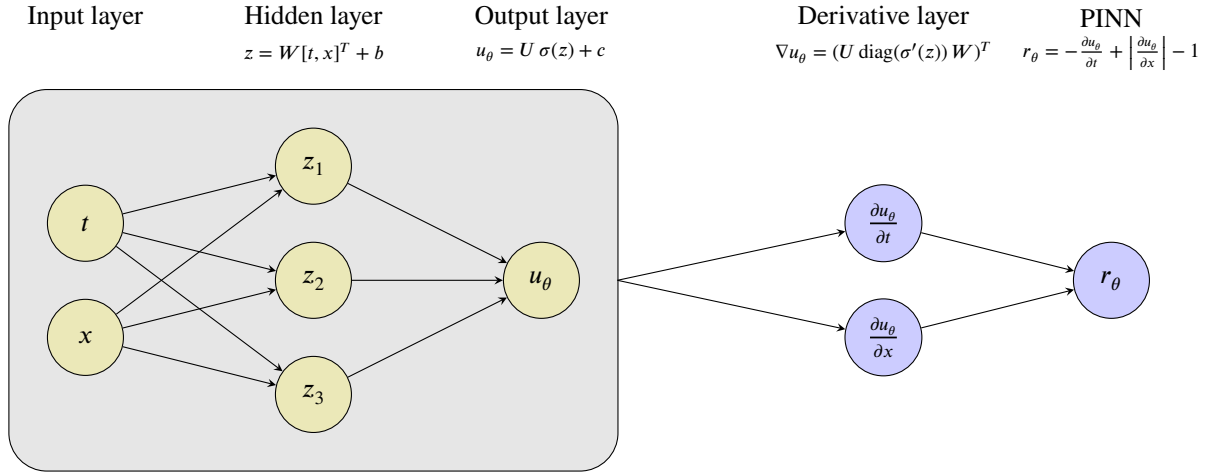

We conclude this section with the explicit derivation of a PINN for a neural network with a single hidden layer.
\begin{example} \label{rem:PINN_Example}
For the one-dimensional eikonal equation~\eqref{eq:Eikonal} the PDE residual is obtained as
\begin{align*}
   r(t,x) := - \partial_t u (t,x) + \abs{\nabla u}(t,x) -1.
\end{align*}
For simplicity we consider a single hidden layer neural network with only three neurons, resulting in the solution approximation 
\begin{align*}
   u_\theta(t,x) = U \sigma\left(W \begin{bmatrix} t\\x \end{bmatrix} + b\right) + c
\end{align*}
with unknown weight matrices $U \in \R^{1 \time 3}, W \in \R^{3 \times 2}$ and bias vectors $b \in \R^3, c \in \R^1$, and an activation function $\sigma:\R \to \R$ acting componentwise on its input.
We further abbreviate the values of the hidden layer by $z = W [t, x]^T + b$.
The chain rule now yields the partial derivatives
\begin{align*}
   \partial_t u_\theta(t,x) = U \, \diag(\sigma'(z))\, W_{:,1}
   \qquad \text{and} \qquad
   \partial_x u_\theta(t,x) = U \, \diag(\sigma'(z))\, W_{:,2}
\end{align*}
where $\diag(\sigma'(z))$ denotes the matrix with diagonal entries $\sigma'(z)$ and $W_{:,j}$ denotes the $j$-th column of the matrix $W$.
This allows us to compute the residual (the actual physics-informed neural network):
\begin{align*}
   r_\theta(t,x) 
   = 
   - U \, \diag(\sigma'(z))\, W_{:,1} + \left| U \, \diag(\sigma'(z))\, W_{:,2} \right| - 1.
\end{align*}
We observe that the residual again possesses the structure of a more complicated neural network mapping $(t,x) \mapsto r(t,x)$.
The neural network employed in this example is illustrated in~\cref{fig:SimpleNN}. 
\end{example}

\subsection{Discrete Time Approach} \label{sec:PINNDisc}

	In contrast to the continuous time approach, the discrete time variant does not incorporate physical information through a set of collocation points, but does so by semi-discretization via Runge-Kutta time-stepping \cite{HairerEtAl1993}.
	Specifically, assuming the solution is known at time $t^n \in [0,T)$, this method assumes the availability of $N_n$ solution data points $X^n:=\{(t^n,x^{n,k}, u^{n,k})\}_{k=1}^{N_n}$ together with boundary data at the domain boundaries.
To continue the solution to $t^{n+1}$, we employ a Runge-Kutta method with $q$ stages
\begin{equation} \label{eq:PINNDiscrete}
    \begin{aligned}
        u^{n+c_i} &= u^n - \Delta t \sum_{j=1}^q a_{ij} \NN [u^{n+c_j}], &&\qquad i=1,\ldots,q,\\
        u^{n+1} &= u^n - \Delta t \sum_{j=1}^q b_j \NN [u^{n+c_j}],
    \end{aligned}
\end{equation}
where $u^{n+c_j} \approx u(t^n + c_j \Delta t, \cdot)$ for $j = 1, \ldots, q$.
Depending on the coefficients $a_{ij}, b_j, c_j$, this represents either an explicit or implicit Runge-Kutta scheme.

While the neural network in the continuous approach approximates the mapping $(t,x) \mapsto u(t,x)$, the discrete-time variant instead approximates $x \mapsto (u^{n+c_1}(x), \ldots, u^{n+c_q}(x), u^{n+1}(x))$, i.e., the solution $u(t,x)$ at the $q+1$ stage values.
Once sufficiently trained, $u(t^{n+1}, x) \approx u^{n+1}(x)$ can be used as the initial data for the next step.
Thus, subsequent steps can proceed analogously.

To be more precise, we establish the link between our data set
$\{(t^n,x^{n,k}, u^{n,k})\}_{k=1}^{N_n}$, the PDE solution at time $t^{n+1}$ and the unknown stages $u^{n+c_i}$, $i=1,\ldots,q$ of the Runge-Kutta scheme~\eqref{eq:PINNDiscrete}, which should hold for all $x \in \domain$, and in particular for all data samples $(x^{n,k}, u^{n,k})$.
This results after a rearrangement of the terms in
\begin{equation*}
    \begin{aligned}
        r^i(x^{n,k}, u^{n,k}) &:= u^{n+c_i}(x^{n,k}) - u^{n,k} + \Delta t \sum_{j=1}^q a_{ij} \NN [u^{n+c_j}](x^{n,k}) \approx 0, &&\qquad i=1,\ldots,q,\\
        r^{q+1}(x^{n,k}, u^{n,k}) &:= u^{n+1}(x^{n,k}) - u^{n,k} + \Delta t \sum_{j=1}^q b_j \NN [u^{n+c_j}](x^{n,k}) \approx 0.
        % u^{n+c_i}(x^{n,k}) &= u^n(x^{n,k}) + \Delta t \sum_{j=1}^q a_{ij} \NN [u^{n+c_j}(x^{n,k})], &&\qquad i=1,\ldots,q,\\
        % u^{n+1}(x^{n,k}) &= u^n(x^{n,k}) + \Delta t \sum_{j=1}^q b_j \NN [u^{n+c_j}(x^{n,k})],
    \end{aligned}
\end{equation*}
These identities are then used to learn the unknown mapping $x \mapsto (u^{n+c_1}(x), \ldots, u^{n+c_q}(x), u^{n+1}(x))$ by minimizing the loss functional, specified here with homogeneous Dirichlet boundary data
\[
    \phi(X^n) 
    := 
    \sum_{k=1}^{N_n}
    % \left(
        \sum_{j=1}^{q+1} \abs{r^j(x^{n,k},u^{n,k})}^2
      %   \abs{u^{n+c_j}(x^{n,k}) - u^{n,k}}^2 
	  % +                 \abs{u^{n+1}(x^{n,k}) - u^{n,k}}^2
    % \right)
    +
    \sum_{i=1}^q
    \big(
        \abs{u^{n+c_i}(-1)}^2
        +\abs{u^{n+c_i}(+1)}^2
    \big)
    +\abs{u^{n+1}(-1)}^2
    +\abs{u^{n+1}(-1)}^2
\]
The numerical experiments presented in \cite{raissi2017physicsI} employ a $\num{500}$-stage Runge-Kutta scheme that advances from initial to final time in a single time step.
The option of using Runge-Kutta methods of extremely high-order rather than small time steps is presented as an attractive feature of this approach, as the task of stage computation for stiff problems requiring implicit integration schemes are passed on to the neural network optimization.
Together with the simplicity of the algorithm and the possibility of choosing large time steps of high order, the numerical results in \cite{raissi2017physicsI} suggest that the method is capable of handling a variety of nonlinearities and boundary conditions.

\subsection{Parameter Identification Setting}

The PINN approach is easily modified to also determine unknown parameters in a general nonlinear partial differential equation.
As an example, consider the PDE
\begin{equation} \label{eq:PINN_ParametricPDE}
    \partial_t u (t,x) + \NN^\lambda[u](t,x) = 0,  \qquad (t,x) \in (0,T] \times \domain,
\end{equation}
with $\NN^\lambda$ a nonlinear partial differential operator depending on a parameter $\lambda \in \R^m$.
Here, we consider only the continuous time framework introduced in \cref{sec:PINNCont}, and refer to~\cite{raissi2017physicsII} for the discrete time variant.

The parameter identification setting as introduced in~\cite{raissi2017physicsII} assumes a set of data $X_d := \{t_i^d,x_i^d, u_i^d\}_{i=1}^{N_d}$, where $u_i^d \approx u(t_i^d, x_i^d)$ are (possibly noisy) observations of the solution of problem~\eqref{eq:PINN_ParametricPDE} in order to identify the unknown parameter $\lambda$.
This training data is then used twofold in a new loss function: in a mean squared misfit term and also in a mean squared residual term:
\[
    \phi(X_d) 
    = 
    \frac{1}{N_d} \sum_{i=1}^{N_d} \abs{u_\theta (t_i^d, x_i^d) - u_i^d}^2 
    + \frac{1}{N_d}  \sum_{i=1}^{N_d} \abs{r_\theta (t_i^d, x_i^d)}^2.
\]

Here, we consider a slightly modified procedure:
In addition to the initial values, boundary and collocation data $X$ introduced in~\cref{sec:PINNCont}, 
we treat the (possibly noisy) observations $X_d$ of the solution of problem~\eqref{eq:PINN_ParametricPDE} in the same way as Dirichlet boundary conditions, which can be enforced via an additional loss function term
\[
	\phi^d(X_d) 
	:= 
	\frac{1}{N_d} \sum_{i=1}^{N_d} \left|u_\theta\left(t_i^d, x_i^d\right) - u_i^d\right|^2,
\]
added to the loss functional~\eqref{eq:PINNLossFunctional}.
% In addition to the initial values, boundary and collocation data $\Xdata$, the framework now requires a set of (possibly noisy) observations of the solution of problem~\eqref{eq:PINN_ParametricPDE} in order to identify the unknown parameters.
% We collect this data in a set $X_d = \{t_i^d,x_i^d, u_i^d\}_{i=1}^{N_d}$.
% This data is treated in the same way as Dirichlet boundary conditions, and can be enforced via an additional loss function term
% \[
%     \phi^d(X_d) 
%     := 
%     \frac{1}{N_d} \sum_{i=1}^{N_d} \left|u_\theta\left(t_i^d, x_i^d\right) - u_i^d\right|^2
% \]
% which is added to the loss functional~\eqref{eq:PINNLossFunctional}.

The unknown parameter $\lambda$ can be learned through training in the same way as the unknown weight matrices $W^\ell$ and bias vectors $b^\ell$ by automatic differentiation of the loss function $\phi$ with respect to $\lambda$.
Indeed, the modifications necessary for including the dependence of the PDE on an unknown parameter require merely a few lines of code, as can be seen in the accompanying Jupyter notebook \texttt{PINN\_Solver.ipynb}.

\begin{figure}[htpb]
    \begin{minipage}{.40\textwidth}
    \centering
        \includegraphics[width=0.98\linewidth]{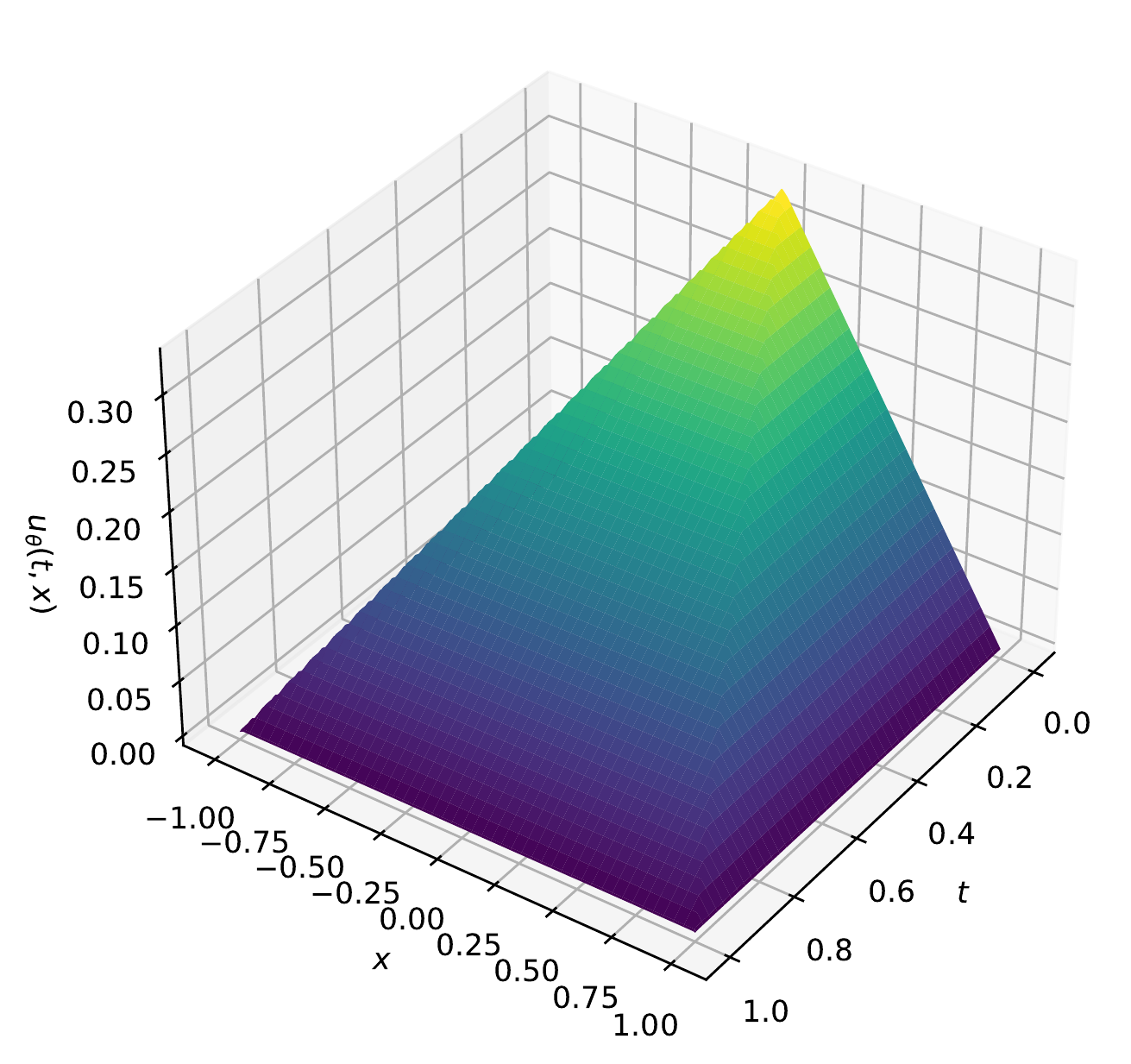}
    \end{minipage}%
    \begin{minipage}{.28\textwidth}
    \centering
        \includegraphics[width=0.98\linewidth]{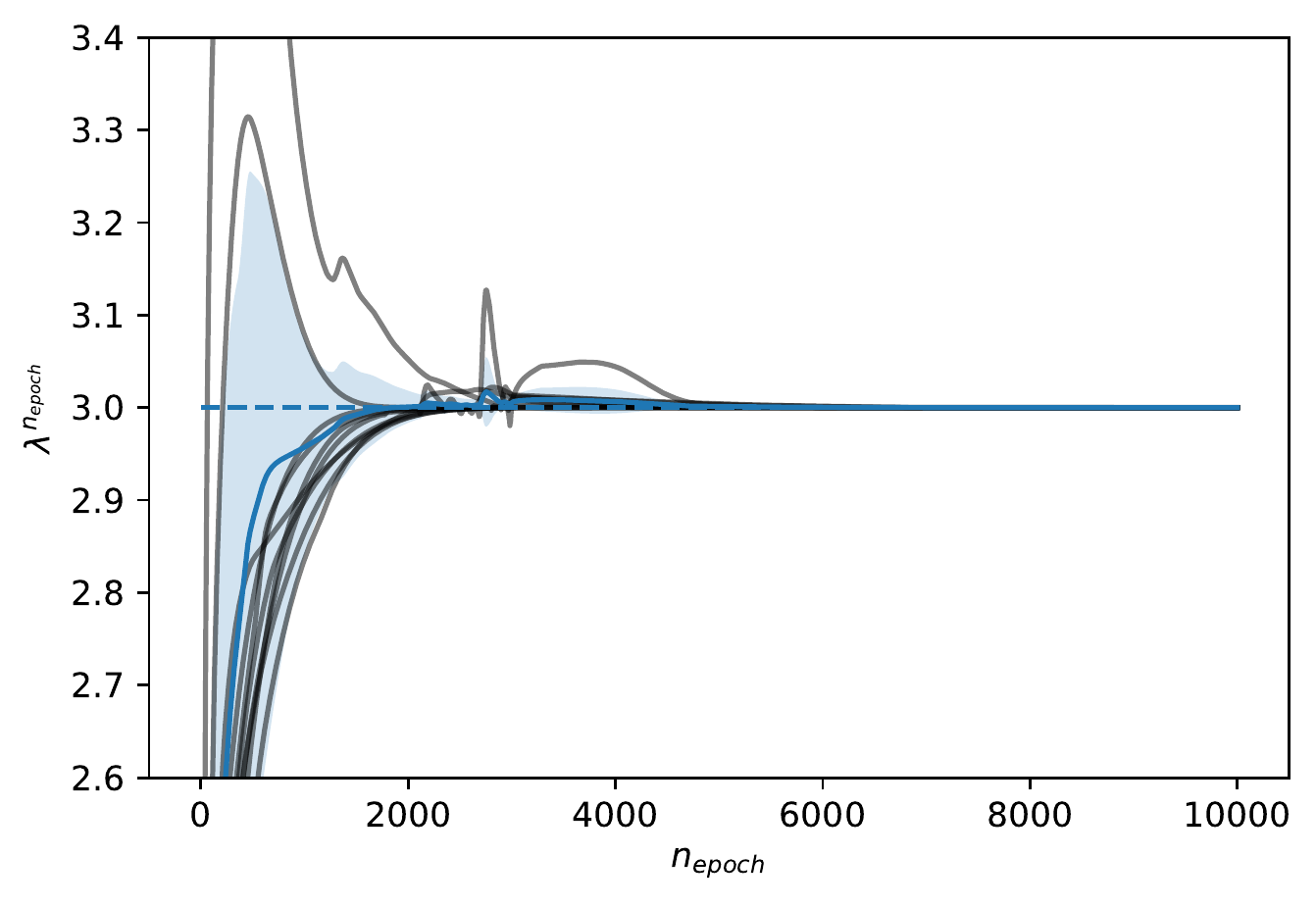}
        \includegraphics[width=0.98\linewidth]{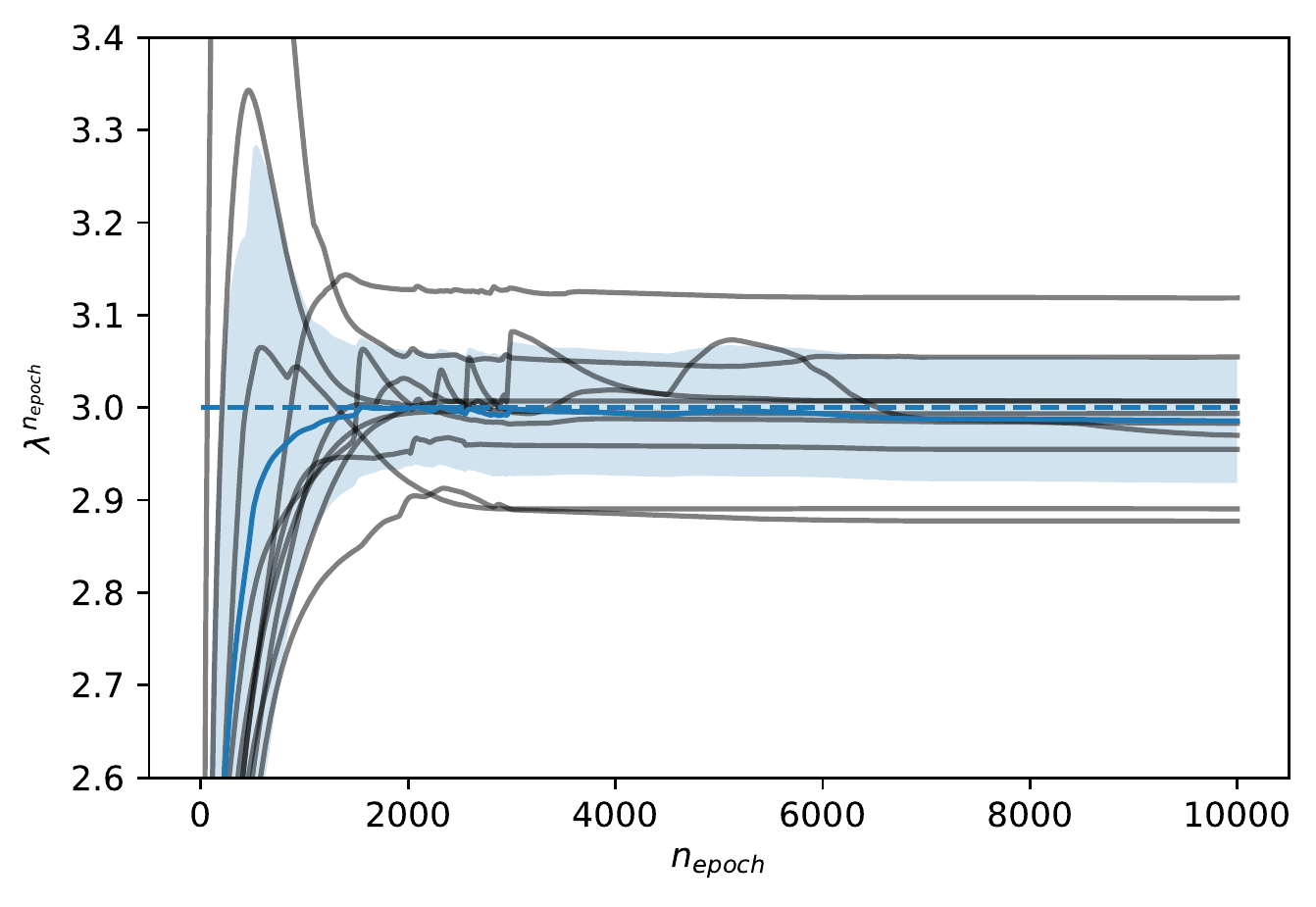}
    \end{minipage}
    \begin{minipage}{.28\textwidth}
    \centering
        \includegraphics[width=0.98\linewidth]{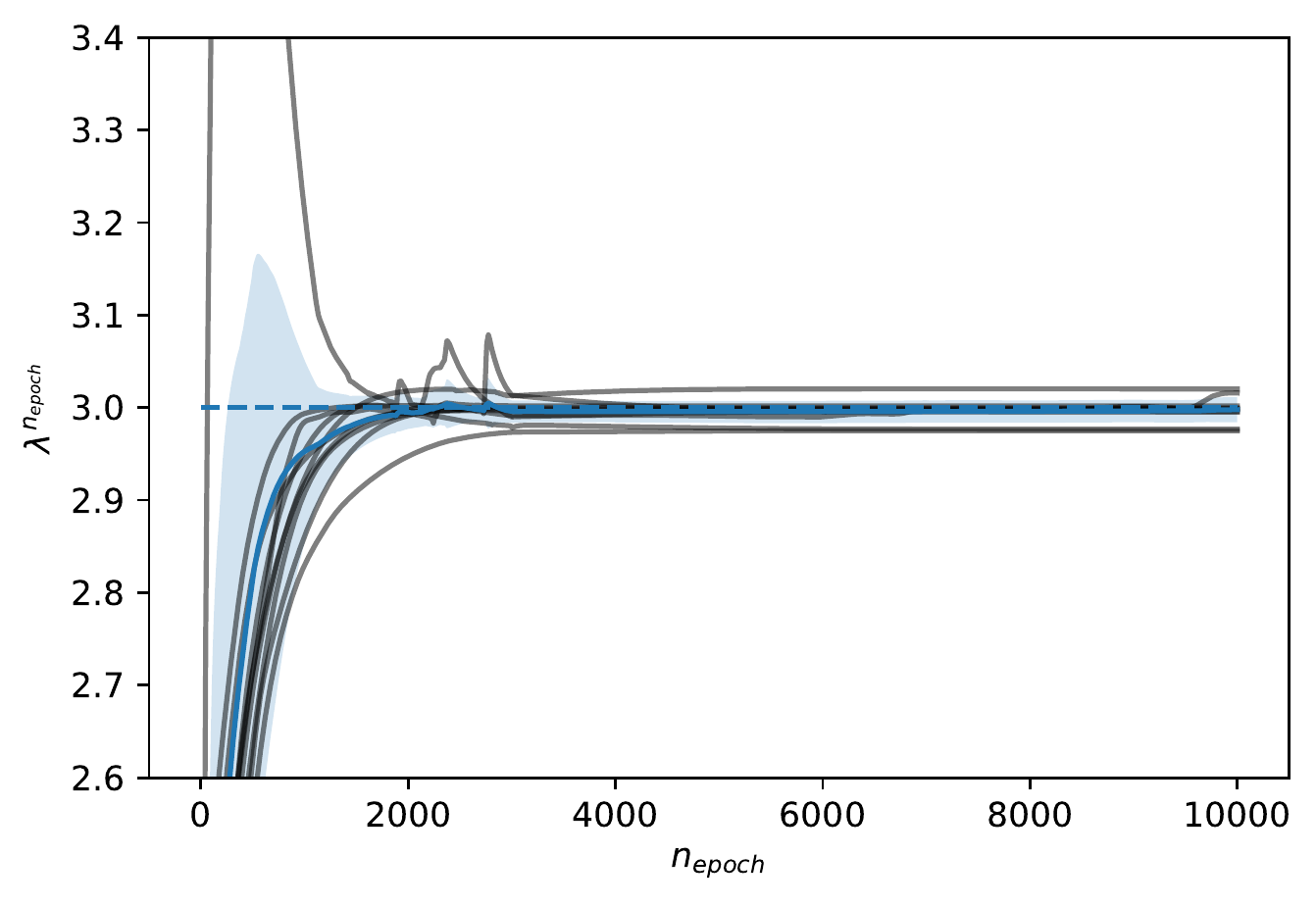}
        \includegraphics[width=0.98\linewidth]{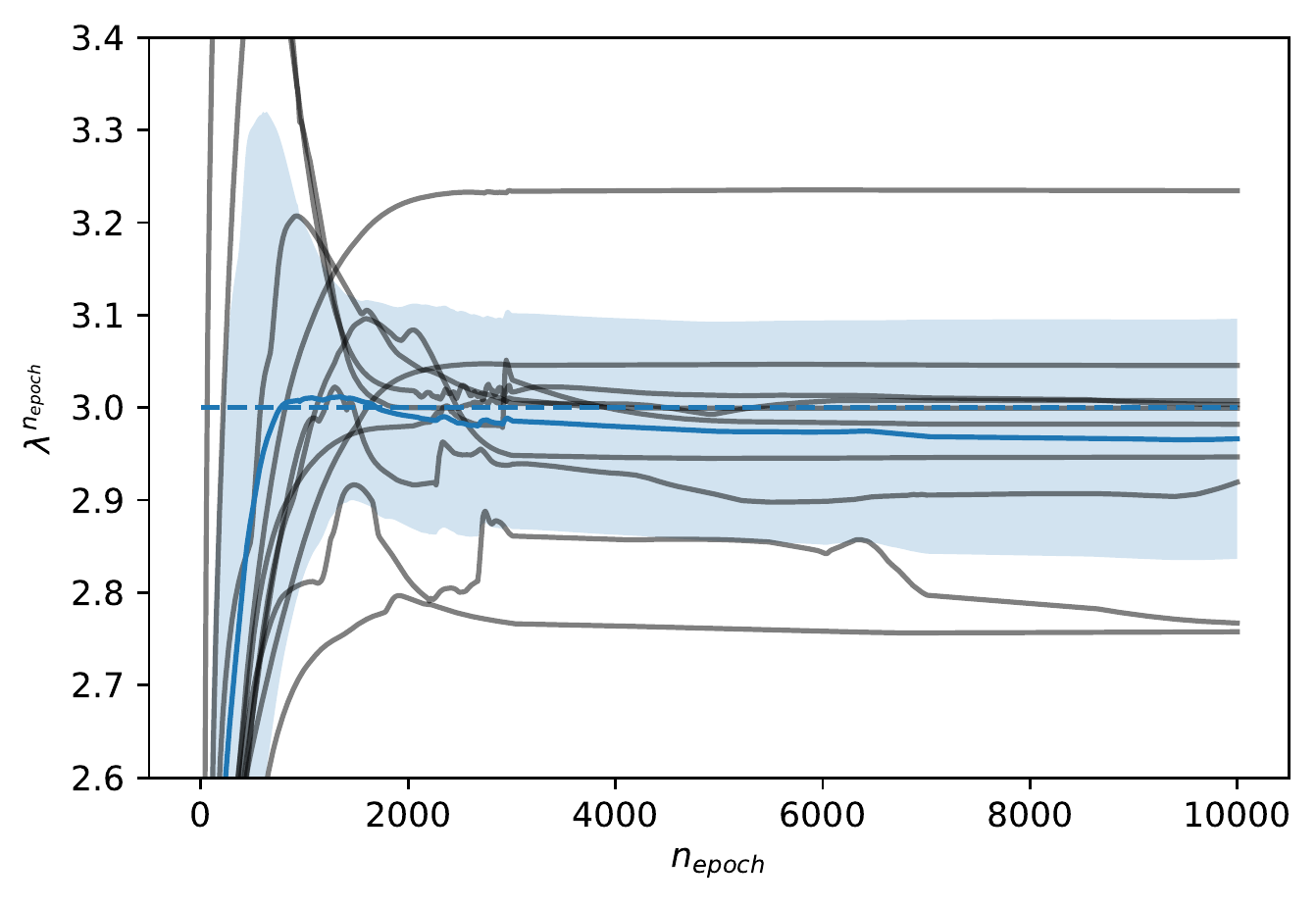}
    \end{minipage}
    \caption{\textbf{Left}: One approximate solution $u_\theta$ of the parametric eikonal equation~\eqref{eq:EikonalParam}.
    \textbf{Right}: Ten evolutions of the estimated parameters $\lambda^{n_\text{epochs}}$ for $n_\text{epochs} = 1, \ldots, \num{10000}$ (gray) with $N_d=500$ noisy measurements $u_i^d = u(t_i^d, x_i^d) + \varepsilon \eta$ with $\eta \sim \pdist{N}(0,1)$ for different noise levels $\varepsilon=0.0, 0.01, 0.05, 0.1$ (from upper left to lower right) together with its mean (solid blue) and one standard deviation around the mean (shaded area).
    The different paths are a result of the random initialization of the parameters in the neural network as well as randomly drawn data $X$ and $X_d$.}%
    \label{fig:Parameteridentification}
\end{figure}

In our example we consider the parametric eikonal equation
\begin{align}
\begin{aligned} \label{eq:EikonalParam}
        -\partial_t u(t,x) + \abs{\nabla u}(t,x) &= \lambda^{-1}
\end{aligned}
\end{align}
with homogeneous final time and boundary conditions and unknown parameter $\lambda > 0$.
Its explicit solution is given by $u^*(t,x) = \lambda^{-1} \, \min\{1-t, 1-\abs{x}\}$.
The numerical results for $\lambda^* = 3$ after $\num{10000}$ training epochs with the Adam optimizer, a piecewise constant learning rate~\eqref{eq:lrEikonal} for a neural network consisting of one hidden layer with $\num{20}$ neurons and leaky ReLU activation function are shown in~\cref{fig:Parameteridentification}.

\subsection{Summary and Extensions} \label{sec:PINNExtension}

Physics-informed neural networks can be used to solve nonlinear partial differential equations.
While the continuous-time approach approximates the PDE solution on a time-space cylinder, the discrete time approach exploits the parabolic structure of the problem to semi-discretize the problem in time in order to evaluate a Runge-Kutta method.
A major advantage of this approach is that it is data-efficient in the sense that it does not require a large number of training samples, which may be difficult to obtain in physical experiments.
Indeed, besides the information on the initial time and spatial boundary, no further knowledge of solution values is required.

In contrast to the method described in \cref{sec:DeepBSDESolver}, the PINN approach is based on a single  neural network to characterize the solution on the entire time-space cylinder $[0,T] \times \overline \domain$.
We note that the focus of the approach does not lie in the solution of high-dimensional problems but rather
in challenging physics features including shocks, convection dominance etc.
Another advantage of this approach is that the value of the loss function can be interpreted as a measure of accuracy of the approximation, and thus can be used as a stopping criterion during training.
We further recall that all derivatives required in the derivation of PINNs~\eqref{eq:PINN} can be computed by the chain rule and evaluated by means of automatic differentiation \cite{baydin2018automatic}.

	A similar physics-constrained approach based on convolutional encoder-decoder neural networks for solving PDEs with random data is developed in  \cite{Zhu2019PhysicsConstrainedDL}.
    Parametrized and locally adaptive activation functions to improve the learning rate in connection with PINNs are explored in \cite{Jagtap2020AdaptiveAF} and \cite{jagtap2020locally}, resp.
	Rigorous estimates on the generalization error of PINNs in the context of inverse problems and data assimilation are given in \cite{mishra2020estimates}.
	XPINNS (eXtended PINNS) are introduced in \cite{JagtapKardiadakis2020} as a generalization of PINNS involving multiple neural networks allowing for parallelization in space and time via domain decomposition, see also \cite{HeinleinKlawonnLanserWeber2020} for a recent review on machine learning approaches in domain decomposition.
	The converse task of learning a nonlinear differential equation from given observations using neural networks is addressed in \cite{RaissiDeepHiddenPhysics}.
	
In addition, PINNs have been applied successfully in a wide range of applications, including
fluid dynamics
\cite{raissi2018hidden,MAO2020112789,lye2020deep,magiera2020constraint,wessels2020neural},
continuum mechanics and elastodynamics \cite{haghighat2020deep,nguyen2020deep,rao2020physics},
inverse problems \cite{meng2020composite,jagtap2020conservative},
fractional advection-diffusion equations
\cite{pang2019fpinns},
stochastic advection-diffusion-reaction equations
\cite{chen2019learning},
stochastic differential equations
\cite{yang2020physics} and
power systems
\cite{misyris2020physics}.
Finally, we mention that Gaussian processes as an alternative to neural networks for approximating complex multivariate functions have also been studied extensively for solving PDEs and inverse problems
\cite{rasmussen2003gaussian,
RaissiEtAl2017TR,%
RaissiPerdikarisKarniadakis2018,%
PangEtAl2019}.
While PINNs have been found to work essentially out of the box in many of these references, as was the case for the examples in \cref{sec:ExampleBurgers}, they may require problem-specific adaptations, particularly when accuracy of efficiency is a consideration. 
An example is a clustering of the interior collocation points to improve the resolution near a shock when solving the Euler equations in \cite{MAO2020112789}.

% Karniadakis et al. Gaussian Processes
% \begin{compactitem}
% \item 
% arXiv preprint \cite{RaissiEtAl2017TR},
% SINUM paper \cite{RaissiPerdikarisKarniadakis2018},
% JCP paper \cite{PangEtAl2019}
% \end{compactitem}

% \jb{Also Gaussian processes and machine learning (Quellen aus Raissi, Karniadakis 2017 TR)
% \cite{rasmussen2003gaussian}
% }

% Recent papers
% \begin{compactitem}
% \item Mao et al.\ 2020 \cite{MAO2020112789} check
% \end{compactitem}

% PINN Vorgänger-Aertikel in GAMM Mitteilungen, aktuell in Revision
% \begin{compactitem}
% \item Axel Klawonn et al.        \cite{HeinleinKlawonnLanserWeber2020}
% \end{compactitem}

% \jb{
% Aus Mao, Karniadakis habe ich noch diese Quellen, alle mit 'check' markierten sind schon verarbeitet, die meisten in dem Satz nach der Aufzählung. Den könntest du einfach noch in die Summary schieben: 
% \cite{rao2020physics} check,
% \cite{lye2020deep} check,
% \cite{magiera2020constraint} check,
% \cite{pang2019fpinns} check,
% \cite{yang2020physics} check,
% \cite{raissi2018hidden} check,
% \cite{chen2019learning} check,
% % \cite{jagtap2020locally} den würde ich mit zu \cite{Jagtap2020AdaptiveAF} (Jagtap2020AdaptiveAF) packen,
% \cite{meng2020composite} check
% }

% \vspace{2em}

\section{Linear PDEs in high Dimensions: the Feynman-Kac Formula} \label{sec:FeynmanKacSolver}

The appeal of the PINN approach of the previous section lies in its simplicity as well as its versatility in applying to a large range of PDE problems. 
The neural network-based approaches presented in this and the next section are aimed at solving PDE problems posed on high-dimensional domains, one of the unsolved problems of numerical analysis.
These problems stem from important applications such as derivative valuation in financial portfolios, the Schrödinger equation in the quantum many-body problem or the Hamilton–Jacobi–Bellman equation in optimal control problems.  
The methods described below are based on the connection between PDEs and stochastic processes, established already in the pioneering work of Bachelier, Einstein, Smoluchowski and Langevin on financial markets, heat diffusion and the kinetic theory of gases (see \cite{Schachermayer2003,Genthon2020} for fascinating accounts) and made explicit in the Feynman-Kac formula \cite{Kac1949}.

In this section and the next, we consider the solution by neural network methods of a class of partial differential equations which arise as the \emph{backward Kolmogorov equation} of stochastic processes known as \emph{It\^{o} diffusions} as proposed in \cite{beck2018solving}.
We begin with linear parabolic second-order partial differential equations in non-divergence form
\begin{equation} \label{eq:KolmogorovPDE}
\begin{aligned}
    \partial_t u(t,x) + \frac{1}{2} \sigma \sigma^T(t,x) : \nabla^2 u(t,x) + \mu(t,x) \cdot \nabla u(t,x) 
    &= 0, 			\quad && (t,x) \in [0,T) \times \R^d,\\
    u(T,x) &= g(x), 	\quad && x \in \R^d,
\end{aligned}
\end{equation}
and subsequently move to more general PDEs.
We consider the pure Cauchy problem, allowing the state variable $x$ to vary throughout $\mathbb R^d$.
Here, $d \in \mathbb{N}$ is the spatial dimension, 
$\nabla u(t,x)$ and $\nabla^2 u(t,x)$ denote the gradient and Hessian of the function $u$, respectively, the colon symbol denotes the Frobenius inner product of $d \times d$ matrices, i.e., $A:B = \sum_{i,j=1}^d a_{ij} \, b_{ij}$, and the dot symbol the Euclidean inner product on $\R^d$.
Let the coefficient functions $\mu\colon[0,T] \times \R^d \to \R^d$ (drift) and $\sigma\colon[0,T] \times\R^d \to \R^{d \times d}$ (diffusion) be globally Lipschitz continuous.
Due to the stochastic process connection, \eqref{eq:KolmogorovPDE} is posed as a \emph{final time problem} with prescribed data at time $t=T$ given by the function $g\colon \R^d \to \R$.
The simple change of variables $t \mapsto T - t$ yields the more familiar initial value form
\begin{equation} \label{eq:KolmogorovInitial}
\begin{aligned}
   \partial_t u(t,x) - \frac{1}{2} \sigma \sigma^T(t,x) : \nabla^2 u(t,x) - \mu(t,x) \cdot \nabla u(t,x) 
   &= 0, 			\quad && (t,x) \in (0,T] \times \R^d,\\
   u(0,x) &= g(x), 	\quad && x \in\R^d.
\end{aligned}
\end{equation}
Equations in non-divergence form like the backward Kolmogorov equation~\eqref{eq:KolmogorovPDE} with  leading term $\sigma \sigma^T(t,x) \colon \nabla^2u(t,x)$ typically arise in the context of stochastic differential equations due to the \Ito\ formula, see \cite{FlemingRishel1975,Ito1944,RevuzYor1999}.
Such problems play a central role in mathematical finance, e.g., in the valuation of complex financial products as well as in stochastic optimal control problems and the solution of second-order Hamilton-Jacobi-Bellman equations \cite{Pham2009,SmearsSuli2014}, where the non-divergence form of the differential operator is again due to the stochastic influence.
Equations of non-divergence type~\eqref{eq:KolmogorovPDE} also arise in the numerical solution of highly nonlinear second-order PDEs that have been linearized, e.g., when applying Newton's method.
Typical examples include the Monge-Ampère equation \cite{BenamouFroeseOberman2010,FengNeilan2009,BrennerNeilan2012}.
Classical and strong solutions of problems in non-divergence form are analyzed in \cite[Ch.\ 6, 9]{GilbargTrudinger2001}.
In contrast to non-divergence PDEs, many problems in applied mathematics arise in divergence form consisting of an operator with leading term 
$\nabla\cdot[\widetilde A(t,x) \nabla u(t,x)]$.
Given sufficient smoothness, each operator in divergence form can be brought into non-divergence form by setting $A(t,x) = \widetilde A(t,x)$ and subtracting the row-wise divergence $\nabla\cdot \widetilde A(t,x)$ from the first-order term.
Even if $\widetilde A$ is smooth, however, this may result in strongly dominating convection in the resulting equation, introducing further challenges.

Following \cite{beck2018solving}, the method reviewed here can be used to construct approximate solutions of a Kolmogorov PDE~\eqref{eq:KolmogorovPDE} or~\eqref{eq:KolmogorovInitial} at a fixed time on some bounded domain of interest $\DD \subset \R^d$.
Similar to the technique reviewed in \cref{sec:PINNs}, a neural network is employed to approximate this solution.
The authors of \cite{beck2018solving} applied their method to a number of examples including the heat equation, the Black-Scholes option pricing equation and others with particular emphasis on the accurate and fast solution in \emph{high dimensions}.
Classical numerical approximation schemes for Kolmogorov partial differential equations are numerous, and include finite difference approximations \cite{brennan1978finite,Kushner1976FD,Kushner1976Survey},
finite element methods \cite{BrennerScott2008,GilbargTrudinger2001,neilan2017convergence,blechschmidt2020error},
numerical schemes based on Monte-Carlo methods \cite{Giles2008,GrahamTalay2013,gobet2016monte,gobet2017adaptive},
as well as approximations based on a discretization of the underlying stochastic differential equations (SDEs) \cite{Higham2011,KloedenPlaten1992}.
Establishing a link of the proposed method to the classical approaches, which might be highly accurate and efficient in up to three dimensions, it shares also similarity to Monte-Carlo methods since it relies on the connection between PDEs and SDEs in the form of the Feynman-Kac theorem and uses a discrete approximation of the SDE associated with equation~\eqref{eq:KolmogorovPDE}.
The reviewed method shares many ideas published in a number of papers, in particular there is a
strong connection to \cite{EHanJentzen2017,HanJentzenE2018} where the Deep BSDE solver, to be presented in detail in \cref{sec:DeepBSDESolver}, is introduced.

In \cite{jentzen2018proof} it is proven that deep neural networks
are able to overcome the curse of dimensionality for
linear backward Kolmogorov PDEs with constant diffusion and nonlinear drift coefficients.
In particular, it is shown that the number of parameters in the neural network grows at most polynomially in both the dimension of the PDE ($d + 1$) and the reciprocal of the desired approximation accuracy.
We note, however, that training a neural network in general is known to be an NP-hard problem, \cite[Sec.\ 20.5]{shalev2014understanding}.

\subsection{The Feynman-Kac Formula} \label{sec:FeynmanKac}

The method reviewed here \cite{beck2018solving} is based on the Feynman-Kac formula for Kolmogorov PDEs which connects the solution of the PDE~\eqref{eq:KolmogorovPDE} and the expectation of a stochastic process.
In order to understand the method fully, we recall the link between PDEs and SDEs formally in this section; for a thorough treatment we refer to \cite{RevuzYor1999,Oksendal2003}.

In a nutshell, the Feynman-Kac theorem states that for every $(t,x) \in [0,T] \times \R^d$ the solution  $u(t,x)$ of the Kolmogorov backward equation~\eqref{eq:KolmogorovPDE} can be expressed as the conditional expectation of a stochastic process $\{X_s\}_{s \in [t,T]}$ starting at $X_t = x$, i.e.,
\begin{equation} \label{eq:ValueU}
    u(t,x) = \E [g(X_T) \given X_t = x].
\end{equation}
Here, $g\colon\R^d \to \R$ is the final time prescribed in \eqref{eq:KolmogorovPDE} and $\E[ \cdot \given X_t=x]$ denotes expectation conditioned on $X_t = x$.
One immediate consequence is that, for all $x \in \R^d$, we have
\begin{equation} \label{eq:ValueBC}
    u(T,x) =  \E [g(X_T) \given X_T = x] = g(x).
\end{equation}
Another implication that can be obtained by the law of iterated conditional expectation is that for all $s \in [t,T]$
\begin{equation} \label{eq:ValueIC}
    u(t,x) = \E [u(s,X_s) \given X_t = x].
\end{equation}

We assume that we are given a filtered probability space $(\Omega, \mathcal{F}, \mathbb{P}; \F)$ equipped with the filtration $\F = \{\FF_t\}_{t \in [0,T]}$ induced by a $d$-dimensional Brownian motion $\{W_t\}_{t \in [0,T]}$.
The stochastic process $\{X_s\}_{s \in [t,T]}$ can be characterized as the solution of the stochastic differential equation (SDE)
\begin{equation} \label{eq:SDE}
    X_s = x + \int_t^s \mu(\tau, X_\tau) \, \d \tau + \int_t^s \sigma(\tau, X_\tau) \, \d W_\tau.
\end{equation}
% See Pham Th. 1.3.15
Assuming Lipschitz conditions on the coefficients $\mu$ and $\sigma$, a pathwise unique strong solution\footnote{Pathwise uniqueness means that if $\{X_s\}_{s\in[t,T]}$ and $\{Y_s\}_{s\in[t,T]}$ are both solutions of~\eqref{eq:SDE}, then $\mathbb P(X_s=Y_s \forall s \in [t,T])=1.$} to~\eqref{eq:SDE} always exists, where $\mu$ and $\sigma$ are the coefficients in~\eqref{eq:KolmogorovPDE}.
Note that the second integral in~\eqref{eq:SDE} is an \Ito\ integral, i.e., a particular type of stochastic integral.
We refer to \cite{RevuzYor1999,Protter2005} for details concerning stochastic analysis and SDEs in general.

Given a strong solution of~\eqref{eq:SDE} $\{X_s\}_{s \in [t,T]}$ and a real-valued function $v \in C^{1,2} ([0,T) \times \R^d; \R) \cap C^0([0,T] \times \R^d; \R)$ applying \Ito's formula \cite{Ito1944,RevuzYor1999}, a generalization of the chain rule for (in generally non-differentiable) stochastic processes,  gives that for any $s \in [t,T]$
\begin{align}    
    v(s,X_s) 
    &= 
    v(t,x) 	+ \int_t^s \partial_t v(\tau, X_\tau) \, \d \tau
    		+ \int_t^s \nabla v(\tau, X_\tau) \cdot \d X_\tau
    		+ \frac{1}{2} \, 
			  \int_t^s \nabla^2 v(\tau, X_\tau):\sigma \sigma^T(\tau, X_\tau) \, \d\tau,
    	\label{eq:uDyn1}
\intertext{which, upon substituting $\d X_\tau$ by its definition~\eqref{eq:SDE}, becomes}
    &= 
    v(t,x) + \int_t^s 
             \left( 
       		 \partial_t v + \frac{1}{2} \nabla^2 v : \sigma \sigma^T + \nabla v \cdot \mu
    		 \right) (\tau, X_\tau) \, \d \tau
           + \int_t^s \nabla v \cdot \sigma(\tau, X_\tau) \, \d W_\tau.
     	\label{eq:uDyn2}
\end{align}
Since this is valid for any $s \in [t,T)$, it holds in particular for $s = t + h$ with $h > 0$, which gives
%together with~\eqref{eq:ValueIC}
\begin{align*}
    v(t+h,X_{t+h}) 
    &= v(t,x) + \int_t^{t+h} 
                \left(
    			\partial_t v + \frac{1}{2} \nabla^2 v : \sigma \sigma^T + \nabla v \cdot \mu
    			\right) (\tau, X_\tau) \, \d \tau
    		  + \int_t^{t+h} \nabla v \cdot \sigma(\tau, X_\tau) \, \d W_\tau.
\end{align*}
Setting $v=u$ given by the expression~\eqref{eq:ValueIC} for $s = t+h$, we obtain
\begin{align*}
    0 &= 
    \E \left[ 
       \int_t^{t+h} 
       \left(
       \partial_t u  + \frac{1}{2} \nabla^2 u : \sigma\sigma^T + \nabla u \cdot \mu
       \right) (\tau, X_\tau) \, \d \tau
       + 
       \int_t^{t+h}  \nabla u \cdot \sigma(\tau, X_\tau) \, \d W_\tau 
    \; \big \vert \; X_t = x \right]\\
    &= 
    \E \left[ 
       \int_t^{t+h} 
       \left( 
       \partial_t u + \frac{1}{2} \nabla^2 u : \sigma\sigma^T + \nabla u \cdot \mu
       \right) (\tau, X_\tau) \, \d \tau
    \; \big \vert \; X_t = x \right],
\end{align*}
where we have used the fact that the stochastic integral is a continuous local martingale and therefore its conditional expectation vanishes.
Dividing by $h > 0$ and taking the limit as $h$ goes to zero yields, by the mean-value theorem,
\begin{align*}
    \partial_t u (t,x)
    + \frac{1}{2} \sigma \sigma^T (t,x) : \nabla^2 u(t,x)
    + \mu(t,x) \cdot \nabla u(t,x)
    = 0
    \qquad \forall (t,x) \in [0,T) \times \R^d,
\end{align*}
confirming that the function given by the Feynman-Kac formula \eqref{eq:ValueIC} solves PDE \eqref{eq:KolmogorovPDE}.

% since they are the result of an application of the \Ito\ formula \cite{Ito1944,RevuzYor1999} applied to a function $u\colon[0,T] \times \R^d \to \R$ that depends on the time $t$ and some stochastic process $\{X_t\}_{t \in [0,T]}$ whose dynamics is given through the \Ito\ process
% $$
% \d X_t = \mu(t, X_t) \, \d t + \sigma(t, X_t) \, \d W_t
% $$
% while $W_t$ is a $d$-dimensional Brownian motion.

% To be more precise, \Ito's formula states that
% \begin{align*}
%     \d u(t, X_t) &= \frac{\partial u}{\partial t} \d t + \nabla u \cdot \d X_t + \frac{1}{2} \nabla^2 u : \d <X_t, X_t>\\
%     &= \frac{\partial u}{\partial t} \d t + \nabla u \cdot 
%     \left(
% \mu(t, X_t) \, \d t + \sigma(t, X_t) \, \d W_t
%     \right)
%     + \frac{1}{2} \nabla^2 u : \sigma \sigma^T(t,X_t) \, \d t,
% \end{align*}
% while $<X_t, X_t>$ denotes the quadratic variation of $X_t$ which is a.s.\ equal to $t$.

\subsection{Methodology} \label{sec:FeynmanKacMethodology}

A number of numerical methods for high-dimensional PDEs have used the Feynman-Kac connection relating PDEs and SDEs in combination with the slow but dimension-independent convergence of Monte Carlo integration 
\cite{%
bouchard2004discrete,%
bender2007forward,%
chassagneux2014linear,%
bender2017primal},
and many more are listed in \cite{BeckEJentzen2019}.
The method from \cite{beck2018solving} reviewed here adds a neural network representation of the PDE solution which is trained in the course of Monte Carlo sampling.
It yields an approximation of the solution $u = u(t,\cdot): \DD \to \R$ of the final time problem \eqref{eq:KolmogorovPDE} restricted to a bounded domain of interest $\DD \subset \R^d$ at a selected time $t \in [0,T]$.
In the following we discuss the methodology in detail for specifically $t=0$.

\subsubsection{Generation of Training Data} \label{sec:FeynmanKacTraining}

Similar to the PINN method discussed in \cref{sec:PINNs}, the method to solve backward Kolmogorov equations does not require any approximate or exact solution values.
Instead, it relies on the generation of a large amount of training data based on the stochastic process connected to the PDE~\eqref{eq:KolmogorovPDE}.

To be more precise, we consider training data $\{(x^i,y^i)\}_{i=1}^{n_\text{data}}$.
Here, the \emph{input} or \emph{independent variable} $x$ is sampled randomly from $X \sim \pdist{U} (\mathcal D)$, which ensures that it covers the domain of interest $\DD$ sufficiently well if sampled many times.
The random \emph{output} (\emph{target variable}) $y$ is defined as a function of $x$ by $Y := g(X_T)$, where $X_T$ is the final value of the stochastic process $\{X_t\}_{t \in[0,T]}$ starting at $X_0=x$ and evolving according to the SDE
\begin{equation}  \label{eq:SDEspec}
	X_t = x + \int_0^t \mu(s, X_s) \, \d s + \int_0^t \sigma(s, X_s) \, \d W_s,
	\qquad 0 \le t \le T.
\end{equation}

We distinguish two cases:

In cases where the distribution of $X_T$ is explicitly known, we can draw sample pairs $(x,y)$ directly.
For example, in the case of a scaled Brownian motion whose dynamics is characterized by $\mu(t,x) \equiv 0$ and $\sigma(t,x) \equiv \sigma$, the solution of~\eqref{eq:SDEspec} is given by
\begin{equation*}
	X_t = x + \sigma \, W_t,
\end{equation*}
where $W_t$ is a path of a standard $d$-dimensional Brownian motion.
Since $W_t \sim \pdist{N}(0, t\, I_{d \times d})$, we may simply draw 
$X \sim \pdist{U}(\mathcal D)$ and set $Y := g(X + \sigma \, \sqrt{T} \, \xi)$, where $\xi \sim \pdist{N}(0, I_{d \times d})$ is a random variable with a $d$-variate standard normal distribution.
Processes for which an explicit distribution is known include Gaussian processes (e.g.\ Brownian motion, Ornstein-Uhlenbeck processes), geometric Brownian motion and Cox-Ingersoll-Ross processes.

When an explicit distribution of $X_t$ at $t=T$ is not available, we may approximate the continuous-time process $\{X_t\}_{t \in [0,T]}$ by generating approximate sample paths using numerical SDE solvers such as the \emph{Euler-Maruyama scheme}
\begin{equation} \label{eq:EulerMaruyama}
	\widetilde X_{n+1} 
	:= 
	\widetilde X_{n}
	+ \mu(t_n, \widetilde X_{n}) \, (t_{n+1} - t_n)
	+ \sigma(t_n, \widetilde X_{n}) \, ( W_{t_{n+1}} -  W_{t_n}),
	\qquad
	\widetilde X_0 := x,
\end{equation}
where $\widetilde X_n \approx  X_{t_n}$ is a discrete-time stochastic process approximating $X_t$ at points $0=t_0 < t_1 < \ldots < t_N = T$ and $x$ is a realization of $X \sim \pdist{U}(\mathcal D)$.
Note that the increment of a Brownian motion $(W_{t_{n+1}} -  W_{t_n}) \sim \pdist{N}(0,(t_{n+1} - t_n) I_{d\times d})$ is normally distributed.
Finally, we set $Y := g(\widetilde X_N)$.
Strong convergence results for the Euler-Maruyama scheme \cite{maruyama1955continuous,KloedenPlaten1992} ensure that $\widetilde X_n \to X_{t_n}$ as $N \to \infty$ and $\sup_n |t_n - t_{n-1}| \to 0$.

\begin{figure}[htpb]
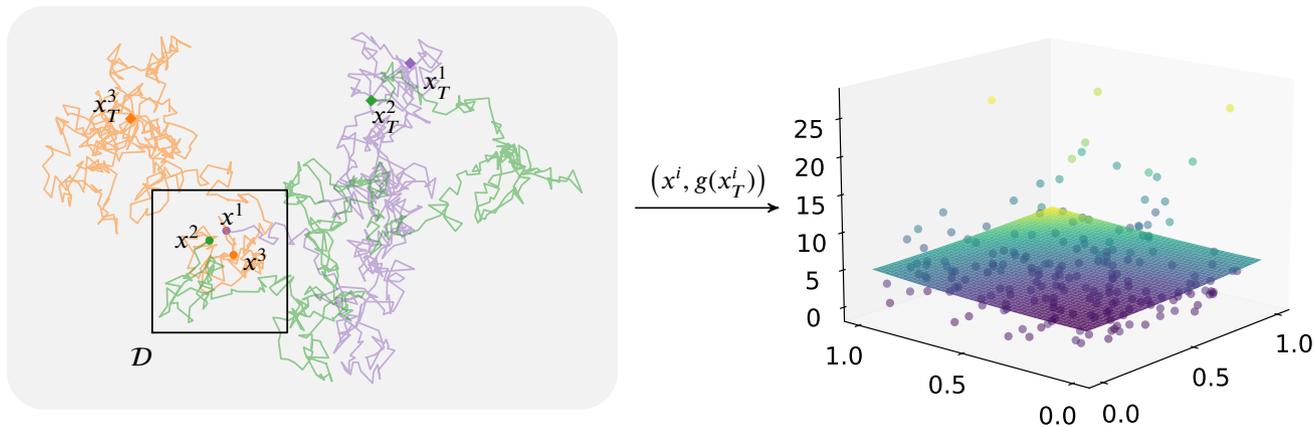

    \centering
\subfile{tikz_FeynmanKac.tex}
    \caption{Illustration of the data generation process described in~\cref{sec:FeynmanKacMethodology}.
    The left panel shows sample paths originating from three different starting values $x^i$ sampled from $\DD$ for $\sigma\equiv\sqrt{2} I_{d \times d}$.
    The right panel shows the exact solution surface $u(0,x)$ along with a number of data pairs $\{(x^i,y^i)\} = \{(x^i, g(x^i_T))\}$ seen to exhibit a large variation around the solution.}%
    \label{fig:IllustrationFeynmanKac}
\end{figure}

	Generating training data via sample paths in this way yields an arbitrary number of easily obtained data pairs $(x^i,y^i)$ with $x^i$ sampled uniformly over $\DD$ and $y^i$ resulting from the final data $g$ evaluated at the final state $X_T$ of a trajectory $\{X_s\}_{s \in [0,T]}$ starting at $X_0=x^i$.
	One has to bear in mind, however, that these individual measurements may vary strongly, in particular for large end times $T$ and diffusion coefficients $\sigma$.
	The training of the neural network $u_\theta: \DD \to \R$ in this way amounts to least squares fitting of $u_\theta$ to a point cloud formed by the data pairs $\{(x^i,y^i)\}_{i=1}^{n_\text{data}}$.
	This is illustrated in~\cref{fig:IllustrationFeynmanKac}.
	The left panel shows sample paths originating from three different starting values $x^i$ sampled from $\DD=[0,1]^2$ for $\sigma\equiv\sqrt{2} I_{d \times d}$.
	Although all processes start within $\DD$, they evolve in $\R^d$ according to the SDE~\eqref{eq:SDEspec} and ultimately leave the domain.
	As a consequence, this method of learning the mapping $u_\theta(0,x)$, $x \in \DD$, does not require the formulation of artificial truncation boundary conditions along $\partial\DD$ as is the case for conventional discretization methods for PDEs on unbounded domains.
	The right panel shows the exact solution surface $u(0,x)$ along with a number of data pairs $\{(x^i,y^i)\} \subset \R^2 \times \R$ seen to exhibit a large variation around the solution.
	Despite the presence of substantial noise in the solution samples, there is no danger of \emph{overfitting} for this method as long as sufficiently many data pairs generated and the training is not restricted to a fixed small number of samples.
	This poses no restriction as the generation of new trajectories and hence solution samples is very inexpensive and allows for an essentially unlimited supply.
	This is particularly true when the distribution of $X_T$ is explicitly known and therefore no numerical path integration is necessary as in the examples given below.

\subsubsection{Neural Network Approximation} \label{sec:FeynmanKacNN}

Similar to the PINN approach discussed in \cref{sec:PINNs}, the unknown solution of the PDE~\eqref{eq:KolmogorovPDE} at a fixed time, here $t=0$, is approximated by a (single) neural network.
We denote this approximation by $u_\theta\colon \DD \to \R$, where $\theta$ collects again all unknown parameters of the network.

The training of the model amounts to a simple regression task.
Given a batch of training data $\{(x^i, y^i)\}_{i=1}^{n_\text{batch}}$, the objective is to minimize the mean squared error
\begin{equation*}
	\frac{1}{n_{\text{batch}}} \sum_{i=1}^{n_{\text{batch}}} |y^i - u_\theta(x^i)|^2,
\end{equation*}
which corresponds from the perspective of the underlying stochastic process to the minimization of
\begin{equation*}
	\mathbb{E} [ |g(X_T) - u_\theta(x)|^2 ]
\end{equation*}
where $X_T$ is the solution of the SDE~\eqref{eq:SDE} starting in $X_0=x$.
This may be viewed as a discrete approximation of a continuous problem, for which it is shown in
\cite[Prop.~2.7]{beck2018solving} 
that, under suitable assumptions, there exists a unique continuous function $u^*:\DD \to \R$ such that
\begin{equation} \label{eq:FeynmanKacObjective}
    \E \left[\abs{g(X_T) - u^*(x)}^2\right] 
    = 
    \inf_{v \in C(\DD;\R)}\E \left[\abs{g(X_T) - v(x)}^2\right].
\end{equation}
Furthermore, it holds for every $x \in \DD$ that $u^*(x) = u(0,x)$.

The network proposed in \cite{beck2018solving}, which is also employed in our numerical tests in \cref{sec:HeatEquation}, has the structure 
\[
   \texttt{Input 
   \lto BN \lto (Dense \lto BN \lto TanH) \lto (Dense \lto BN \lto TanH) 
   \lto Dense \lto BN
   \lto Output}
\]
where the notation is as follows:
\begin{itemize}
\item 
	\texttt{BN} indicates a  \emph{batch normalization} step \cite{ioffe2015batch}, which is a technique of normalizing each training mini-batch within the network architecture to make the model less sensitive in terms of proper weight initialization and allows for larger step sizes and faster training.
        This is effected by additional parameters that scale and shift the neurons that enter the \texttt{BN} layer componentwise.
        These parameters are learned in the same way as all unknown parameters in the neural network, e.g., by a mini-batch gradient descent type algorithm.
\item
	\texttt{Dense} indicates a fully connected layer \emph{without} bias term, i.e., a matrix-vector product with a learnable weight matrix.
	Due to the subsequent shifting during the \texttt{BN} layer, a bias term can be omitted since its effect would be cancelled.
\item
	\texttt{TanH} indicates the application of the componentwise hyperbolic tangent activation function.
% \begin{equation*}
%     \tanh(x) = \frac{\e^x - \e^{-x}}{\e^x + \e^{-x}}.
% \end{equation*}
\end{itemize}

The network is trained with the Adam optimizer \cite{kingma2014adam}, a variant of the stochastic gradient descent method based on an adaptive estimation of first-order and second-order moments to improve the speed of convergence.
An explanatory walkthrough of the implementation of the complete algorithm is given in the accompanying Jupyter notebook \texttt{Feynman\_Kac\_Solver.ipynb}.

\subsubsection{Example: Heat equation} \label{sec:HeatEquation}

In this section, we want to solve the heat equation in $d$ dimensions by means of the solver proposed in~\cite{beck2018solving} and consider the initial value problem
\begin{equation} \label{eq:HeatEquation}
    \begin{aligned}
        \partial_t u(t,x) &= \Delta u (t,x) \quad & (t,x) \in (0,T] \times \R^d\\
        u(0,x) &= \|x\|^2 \quad & x \in \R^d,
    \end{aligned}
\end{equation}
where $\Delta u = \sum_{i=1}^d \partial^2 u / \partial x^2_i$ denotes the Laplacian of $u$.
One can easily verify that the solution is given by
\begin{equation*}
    u(t,x) = \|x\|^2 + 2 \, t \, d.
\end{equation*}

% \jb{Var.\ A:
% During our numerical tests we observed that the convergence behaviour is highly sensitive.
% At a first glance, it seems tempting to choose a larger initial step size which results in a faster decay of the errors in the beginning, however, the decay stagnates earlier and the final results are inferior.}
% \jb{Var\ B:
We tested two different step size strategies: a decaying piecewise constant learning rate with step sizes $\delta(n) = 10^{-3} \, \textbf{1}_{\{n \le \num{250000}\}} + 10^{-4} \, \textbf{1}_{\{\num{250000} < n \le \num{500000}\}} + 10^{-5} \, \textbf{1}_{\{\num{500000} < n\}}$ as was employed in \cite{beck2018solving} and an exponentially decaying rate with step sizes $\delta(n) = 0.1 \cdot 10^{-n/\num{100000}}$.
The remainder of the parameters are chosen as in~\cite{beck2018solving}.
We fixed the number of neurons in the two hidden layers to $\num{200}$ independent of the dimension.
\Cref{fig:HeatErrors} shows the evolution of the absolute and relative approximation errors\footnote{All errors shown in the plots are approximated by Monte-Carlo estimation with one million samples.} on $\DD$ for the $100$-dimensional heat equation.

\pgfplotstableread[col sep = comma]{heat_dim_100_orig.csv}\TabOrig%
\pgfplotstableread[col sep = comma]{heat_dim_100_exp.csv}\TabExp%

\begin{figure}[htpb]
    \centering
    \begin{subfigure}{0.30\textwidth}
		\centering
		\resizebox{\textwidth}{!}{
			\begin{tikzpicture}
				\begin{semilogyaxis}[%
                        enlargelimits=false,
                    xlabel={$n_\text{epoch}$},
                    label style={font=\footnotesize}]
                    \addplot+[L1relerror, solid] table[x=Iter,y=L1_rel] {\TabOrig};
                    \addplot+[L2relerror, solid] table[x=Iter,y=L2_rel] {\TabOrig};
                    \addplot+[Linfrelerror, solid] table[x=Iter,y=Linf_rel] {\TabOrig};
                    \addlegendentry{$L^1_{\text{rel}}$};
					\addlegendentry{$L^2_{\text{rel}}$};
					\addlegendentry{$L^\infty_{\text{rel}}$};
                    \addplot+[L1relerror, dashed] table[x=Iter,y=L1_rel] {\TabExp};
                    \addplot+[L2relerror, dashed] table[x=Iter,y=L2_rel] {\TabExp};
                    \addplot+[Linfrelerror, dashed] table[x=Iter,y=Linf_rel] {\TabExp};
				\end{semilogyaxis}
			\end{tikzpicture}
		}
    \end{subfigure}
    \hfill
    \begin{subfigure}{0.39\textwidth}
        \includegraphics[width=0.98\textwidth]{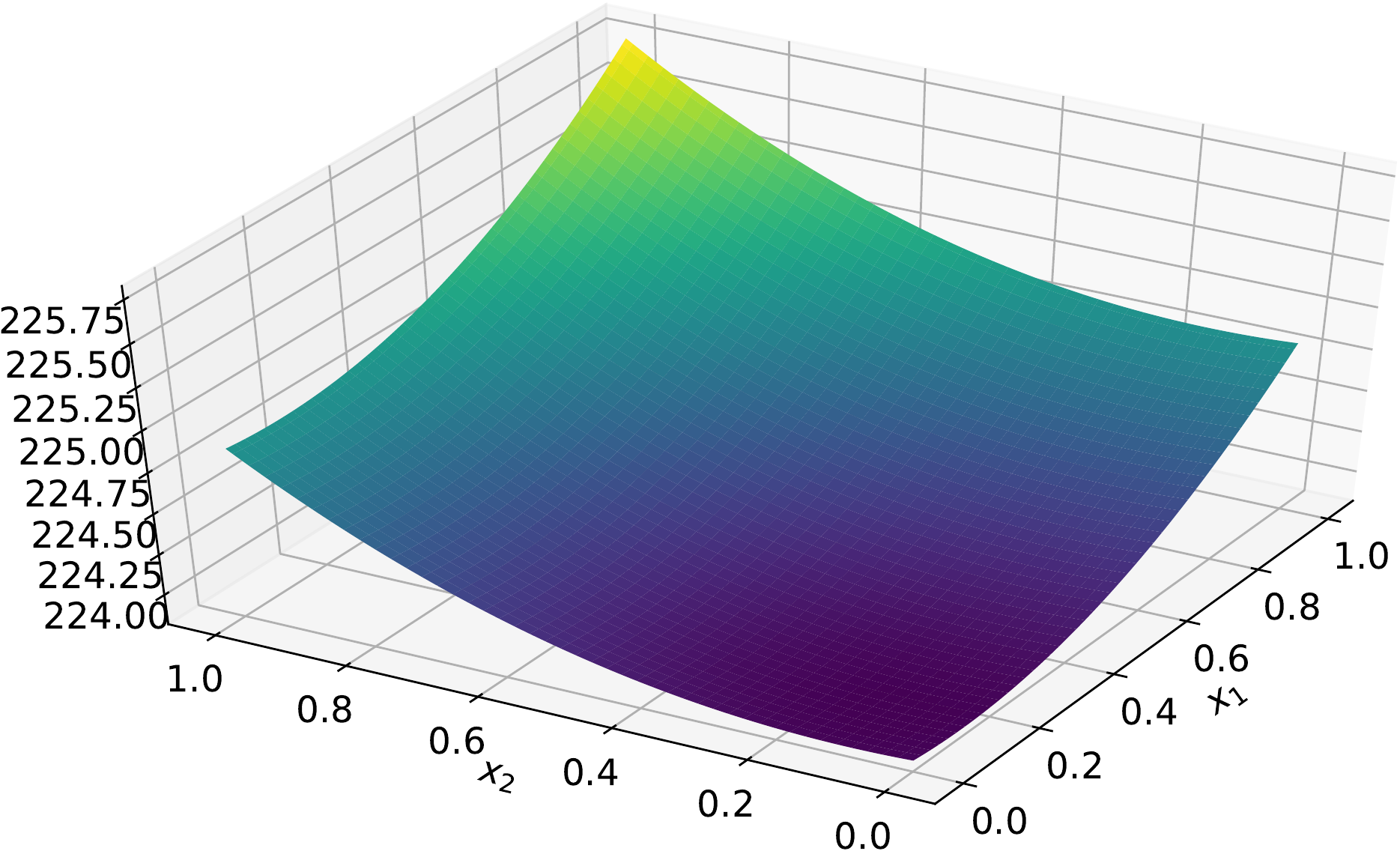}
    \end{subfigure}
    \hfill
    \begin{subfigure}{0.30\textwidth}
		\centering
		\resizebox{\textwidth}{!}{
			\begin{tikzpicture}
				\begin{semilogyaxis}[%
                        enlargelimits=false,
                    xlabel={$n_\text{epoch}$},
                    label style={font=\footnotesize}]
                    \addplot+[L1abserror, solid] table[x=Iter,y=L1_abs] {\TabOrig};
                    \addplot+[L2abserror, solid] table[x=Iter,y=L2_abs] {\TabOrig};
                    \addplot+[Linfabserror, solid] table[x=Iter,y=Linf_abs] {\TabOrig};
					\addlegendentry{$L^1_{\text{abs}}$};
					\addlegendentry{$L^2_{\text{abs}}$};
					\addlegendentry{$L^\infty_{\text{abs}}$};
                    \addplot+[L1abserror, dashed] table[x=Iter,y=L1_abs] {\TabExp};
                    \addplot+[L2abserror, dashed] table[x=Iter,y=L2_abs] {\TabExp};
                    \addplot+[Linfabserror, dashed] table[x=Iter,y=Linf_abs] {\TabExp};
				\end{semilogyaxis}
			\end{tikzpicture}
		}
    \end{subfigure}
    \caption{Evolution of relative  (\textbf{left}) and absolute  (\textbf{right}) errors for a decaying piecewise constant  learning rate (solid) and an exponentially decaying rate (dashed) for the $\num{100}$-dimensional heat equation~\eqref{eq:HeatEquation}, estimated by means of the Monte-Carlo method in order to approximate the integrals in dimension $\num{100}$ with one million samples. \textbf{Center}: Two-dimensional slice through the approximate solution $(x_1,x_2) \mapsto u_\theta(x_1,x_2,0.5,\ldots,0.5)$.}
    \label{fig:HeatErrors}
\end{figure}
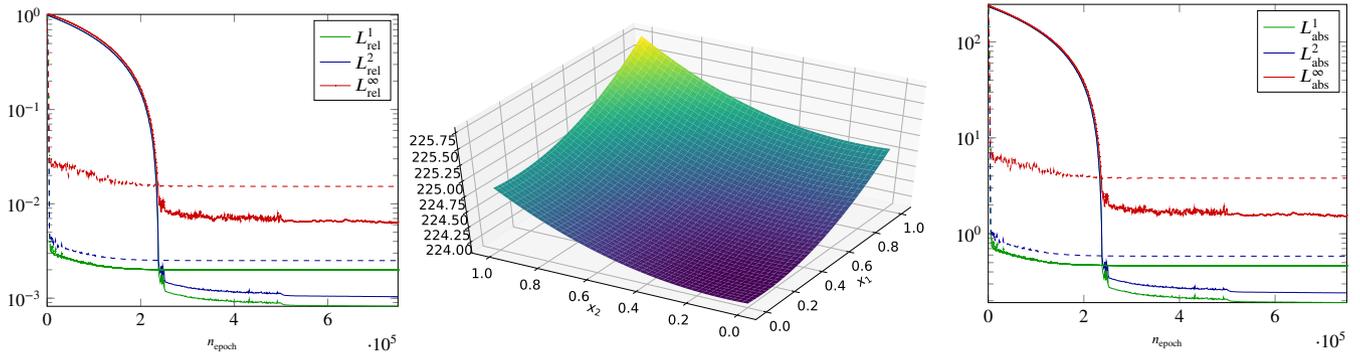
	% \bigskip

%     \centering
%     \hypersetup{hidelinks}
%     % \ref{exp3degreeStudyLegend} % Here comes the legend
%     \label{fig:experiment3_degree_study}
% \end{figure}

% \begin{figure}[htpb]
%     \centering
%     \caption{\textbf{Left}: Shown are the evolutions of relative errors (solid) and absolute errors (dashed) for the $\num{100}$-dimensional heat equation~\eqref{eq:HeatEquation}, estimated by means of the Monte-Carlo method in order to approximate the integrals in dimension $\num{100}$ with one million samples. \textbf{Right}: Shown is a two dimensional slice through the approximate solution $(x_1,x_2) \mapsto u_\theta(x_1,x_2,0.5,\ldots,0.5)$.}%
%     \label{fig:HeatErrors}
% \end{figure}

	In our numerical experiments with the heat equation~\eqref{eq:HeatEquation} we observed that the quality of the final approximation depends heavily on the chosen learning rate, i.e., the step sizes used in the gradient method.
	A comparison between the evolutions of the relative and absolute errors for the two aforementioned learning rate strategies is displayed in \cref{fig:HeatErrors}, together with a two-dimensional slice through the 100-dimensional solution.
	Together with \cref{table1}, this indicates that it seems to be better to stay conservative and take smaller steps from the beginning on.
Shown are errors for the two step size scenarios at $\nepoch = \num{100000}$ and $\nepoch = \num{750000}$.
Although the errors decrease faster in the beginning for the exponentially decaying step sizes that start with larger steps, the errors seems to saturate at a higher level.
	This might be due to the algorithm settling into some local minimum.
For the decaying piecewise constant learning rate, \cref{fig:HeatErrors} shows two distinct phases of error decay: While the first phase until approximately epoch number $\num{250000}$ is characterized by an accelerating decay of the errors probably due to mainly shifting the solution slowly towards the image range ($200$ to $300$) of the solution, the second phase decays at a much slower rate which might correspond to the reduction rate of the Monte Carlo error.
The exponentially decaying learning rate decays much faster in the beginning but settles at a higher absolute and relative error.

\pgfplotstableset{
    columns/id/.style={
        column name=Experiment,
        string type},
    multicolumn names, % allows to have multicolumn names
    col sep=comma, % the seperator in our .csv file
    sci zerofill,
    columns/Niter/.style={
        column name=Dim, % name of first column
        int detect},  % use siunitx for formatting
    columns/L1rel/.style={
        column name=$L_{\text{rel}}^1(\mathcal{D})$,
        sci},
    columns/L2rel/.style={
        column name=$L_{\text{rel}}^2(\mathcal{D})$,
        sci},
    columns/Linfrel/.style={
        column name=$L_{\text{rel}}^\infty(\mathcal{D})$,
        sci},
    columns/L1abs/.style={
        column name=$L_{\text{abs}}^1(\mathcal{D})$,
        sci},
    columns/L2abs/.style={
        column name=$L_{\text{abs}}^2(\mathcal{D})$,
        sci},
    columns/Linfabs/.style={
        column name=$L_{\text{abs}}^\infty(\mathcal{D})$,
        sci},
    columns/Time/.style={
        column name=Time,
        column type={S},string type},
    % ,
    %     sci e,
    %     sci zerofill},
    % columns/NdofsMixed/.style={
    %     column name=$\dim(V_h + W_h^{2 \times 2})$,
    %     int detect},
    % columns/L2_error/.style={
    %     column name=$\norm{u-u_h}_{L^2(\Omega)}$,
    %     sci,
    %     sci e,
    %     sci zerofill},
    % columns/L2_eoc/.style={
    %     string replace={0}{},
    %     column name=$\text{EOC}_{L^2}$,
    %     precision=2,
    %     fixed,
    %     fixed zerofill,
    %     postproc cell content/.append style={/pgfplots/table/@cell content/.add={(}{)},},
    %     string replace={()}{},
    % },
    % columns/H1_error/.style={
    %     column name=$\norm{u-u_h}_{H^1_0(\Omega)}$,
    %     sci,
    %     sci e,
    %     sci zerofill},
    % columns/H1_eoc/.style={
    %     string replace={0}{},
    %     column name=$\text{EOC}_{H^1_0}$,
    %     sci,
    %     sci zerofill,
    %     precision=2, fixed},
    % columns/H2_error/.style={
    %     column name=$\norm{u-u_h}_{H^2_0(\Omega)}$,
    %     sci,
    %     sci e,
    %     sci zerofill},
    % columns/H2_eoc/.style={
    %     string replace={0}{},
    %     column name=$\text{EOC}_{H^2_0}$,
    %     precision=2,
    %     fixed,
    %     fixed zerofill},
    % columns/H2h_error/.style={
    %     column name=$\norm{u-u_h}_{H^2_h(\Omega)}$,
    %     sci,
    %     sci e,
    %     sci zerofill},
    % columns/H2h_eoc/.style={
    %     string replace={0}{},
    %     column name=$\text{EOC}_{H^2_h}$,
    %     sci,
    %     sci zerofill,
    %     precision=2, fixed},
    % columns/gmres_iter/.style={
    %     column name=$N_\text{gmres}$,
    %     column type={S},string type},
    every head row/.style={
        before row={\toprule}, % have a rule at top
        after row={\midrule}, % rule under units
    },
    every nth row={3}{before row=\midrule},
    every last row/.style={after row=\bottomrule}, % rule at bottom
    }
    \begin{table}[htb]
        \begin{center}
            \begin{filecontents*}{./table_Heat.csv}
id,Dim,Iter,Loss,L1rel,L2rel,Linfrel,L1abs,L2abs,Linfabs,Time,Stepsize
,10,250000,110.76896,0.0015347758,0.0019664331,0.030855268,0.035983644,0.04646216,0.857502,896.447301864624,0.0003162205
Exp.\ decay,50,250000,529.32776,0.0021366964,0.0026950894,0.026115084,0.24923484,0.31442472,3.3424454,940.5341627597809,0.0003162205
$n_\text{epoch}=\num{250000}$,100,250000,1051.7157,0.0019680473,0.0024691646,0.013513848,0.45915475,0.5761033,3.3485107,1038.7952525615692,0.0003162205
,10,750000,104.781296,0.0014868401,0.0019059435,0.030654315,0.034865737,0.045052182,0.85191727,2692.6939449310303,3.1622054e-09
Exp.\ decay,50,750000,548.8055,0.0021179777,0.002671101,0.02607604,0.24703892,0.31160578,3.3374481,2830.4973883628845,3.1622054e-09
$n_\text{epoch}=\num{750000}$,100,750000,1054.8625,0.001955854,0.0024544718,0.013247941,0.45631963,0.5726988,3.2826233,3129.64191865921,3.1622054e-09
,10,250000,104.77309,0.002402247,0.003006398,0.013023739,0.05598822,0.0701832,0.3448429,1122.7159886360168,0.001
Piecewise decay,50,250000,532.4577,0.0015896752,0.001998788,0.011028128,0.18540815,0.23306017,1.2491074,1153.8490903377533,0.001
$n_\text{epoch}=\num{250000}$,100,250000,1074.5591,0.0014398925,0.0018056053,0.009659108,0.3357871,0.42095923,2.2595215,1203.9419221878052,0.001
,10,750000,106.98634,0.00059044466,0.0007427701,0.004288846,0.01374449,0.01726429,0.091947556,3405.134256362915,1e-05
Piecewise decay,50,750000,534.0857,0.00077584747,0.0009869041,0.007151018,0.090430625,0.114959754,0.8661194,3493.8051540851593,1e-05
$n_\text{epoch}=\num{750000}$,100,750000,1042.1028,0.0008204926,0.0010374382,0.006244784,0.19133218,0.24183008,1.4659424,3659.943394422531,1e-05
\end{filecontents*}
\pgfplotstabletypeset[
                % columns={{Iter},{L1_rel},{L2_rel},{Linf_rel},{L1_abs},{L2_abs},{Linf_abs},{Time}},
                % columns={{Dim},{L1_rel},{L2_rel},{Linf_rel},{L1_abs},{L2_abs},{Linf_abs},{Time}},
                columns={{id},{Dim},{L1rel},{L2rel},{Linfrel},{L1abs},{L2abs},{Linfabs},{Time}},
            ]{./table_Heat.csv} % filename/path to file
            % Iter,Loss,L1_rel,L2_rel,Linf_rel,L1_abs,L2_abs,Linf_abs,Time,Stepsize
        \end{center}
            \caption{Absolute and relative approximation errors for the $d$-dimensional heat equation~\eqref{eq:HeatEquation}.}
            \label{table1}
    \end{table}

We also observe that it seems to be difficult to improve the achievable relative and absolute errors, see \cref{table:Comparison}.

\pgfplotstableset{
    columns/id/.style={
        column name=Experiment,
        string type},
    multicolumn names, % allows to have multicolumn names
    col sep=semicolon, % the seperator in our .csv file
    sci zerofill,
    columns/Niter/.style={
        column name=Dim, % name of first column
        int detect},  % use siunitx for formatting
    columns/L1rel/.style={
        column name=$L_{\text{rel}}^1(\mathcal{D})$,
        sci},
    columns/L2rel/.style={
        column name=$L_{\text{rel}}^2(\mathcal{D})$,
        sci},
    columns/Linfrel/.style={
        column name=$L_{\text{rel}}^\infty(\mathcal{D})$,
        sci},
    columns/L1abs/.style={
        column name=$L_{\text{abs}}^1(\mathcal{D})$,
        sci},
    columns/L2abs/.style={
        column name=$L_{\text{abs}}^2(\mathcal{D})$,
        sci},
    columns/Linfabs/.style={
        column name=$L_{\text{abs}}^\infty(\mathcal{D})$,
        sci},
    columns/Time/.style={
        column name=Time,
        column type={S},string type},
    every head row/.style={
        before row={\toprule}, % have a rule at top
        after row={\midrule}, % rule under units
    },
    every nth row={3}{before row=\midrule},
    every last row/.style={after row=\bottomrule}, % rule at bottom
    }
    \begin{table}[htb]
        \begin{center}
            \begin{filecontents*}{./table_100_comp.csv}
id;Dim;Iter;Loss;L1rel;L2rel;Linfrel;L1abs;L2abs;Linfabs;Time;Stepsize
$n_\text{layers}=2,n_\text{neuron}=200$;100;750000;1042.1028;0.0008204926;0.0010374382;0.006244784;0.19133218;0.24183008;1.4659424;3659.943394422531;1e-05
$n_\text{layers}=3,n_\text{neuron}=300$;100;750000;1093.5913;0.0007754837;0.0009786362;0.0066411737;0.18083414;0.22809455;1.5127106;4704.552524805069;1e-05
$n_\text{layers}=4,n_\text{neuron}=400$;100;750000;1076.6865;0.0007237474;0.00091490196;0.0054260325;0.168773;0.21325718;1.3014374;8483.335824251175;1e-05
\end{filecontents*}
\pgfplotstabletypeset[
                columns={{id},{L1rel},{L2rel},{Linfrel},{L1abs},{L2abs},{Linfabs},{Time}},
            ]{./table_100_comp.csv} % filename/path to file
        \end{center}
            \caption{Absolute and relative errors of the 100-dimensional heat equation with decaying piecewise constant learning rate for three different neural network architectures after $\num{750000}$ training epochs.}
            \label{table:Comparison}
    \end{table}

In this example no SDE time-stepping is necessary, as the end of the sample paths $X_T$ can be drawn directly. 
In particular, this incurs no discretization error.

A general recomendation on how to select the neural network architecture and parameter selection could be part of further research.
This however, is a problem prevalent in many fields of scientific machine learning, see \cite{petersen2018optimal} for a discussion on selecting deep ReLU network architectures.
Nevertheless, one has to bear in mind that problems in such a high spatial dimension have been considered absolutely infeasible for a long time in terms of numerical approximations.
In particular, for problems in financial mathematics where derivatives, e.g., options, often depend on a basket of more than 100~underlying risky assets (which determine the spatial dimension of the pricing PDE), the importance of having a feasible algorithm can not be denied.
Note that the accompanying code includes as a second example an option pricing problem.

\subsection{Linear Parabolic PDEs in General Form} \label{sec:FeynmanKacGeneral}

The Feynman-Kac formula may be extended to the full class class of linear parabolic equations, see~\cite[Ch.~5 Theorem~7.6]{karatzas2014brownian}.
Specifically, adding a zeroth order term with non-negative potential $r \colon [0,T] \times \R^d \to [0,\infty)$ as well as a source term $f \colon [0,T] \times \R^d \to \R$, the final time problem \eqref{eq:KolmogorovPDE} becomes
\begin{equation} \label{eq:LinearParabolicPDE}
\begin{aligned}
    \partial_t u(t,x) 
    + \frac{1}{2} \sigma \sigma^T(t,x) : \nabla^2 u(t,x) 
    + \mu(t,x) \cdot \nabla u(t,x)  
    - r(t,x) \, u(t,x)
    + f(t,x)
    &= 0, 			\quad && (t,x) \in [0,T) \times \R^d,\\
    u(T,x) &= g(x), 	\quad && x \in \R^d.
\end{aligned}
\end{equation}
A sufficiently smooth solution of \eqref{eq:LinearParabolicPDE} admits the Feynman-Kac representation 
\begin{equation} \label{eq:FeynmanKacGeneral}
    u(t,x) 
    = 
    \E \left[
        \int_t^T \e^{- \int_t^\tau r(\nu,X_\nu) \, \d \nu} \, f(\tau,X_\tau) \, \d \tau
        + 
        \e^{- \int_t^T r(\nu,X_\nu) \, \d \nu} \, g(X_T) \given X_t = x
        \right]
        \qquad
        \forall (t,x) \in [0,T] \times \R^d,
\end{equation}
which simplifies to \eqref{eq:ValueU} for $f \equiv 0$ and $r \equiv 0$.

Algorithmically, this can be considered within the same framework as discussed in~\cref{sec:FeynmanKacMethodology}.
% Algorithmically, this can be considered within the same framework as discussed in~\cref{sec:FeynmanKacMethodology}.
In particular, it does not change the generation of samples of the stochastic process $\{X_s\}_{s \in [0,T]}$.
In the case of a discrete approximation $\{\widetilde X_n\}_{n=0}^N$ generated by the Euler-Maruyama scheme~\eqref{eq:EulerMaruyama},
a simple approximation of the corresponding output variable $Y$ can be given by
\begin{equation} \label{eq:YProcessFull}
    Y 
    = 
    \sum_{n=0}^{N-1} 
        \widetilde R_{n} \, f(t_n,\widetilde X_n) \, (t_{n+1}-t_n) 
      + \widetilde R_{N} \, g(\widetilde X_N)
\end{equation}
with
\[
    \widetilde R_n 
    := \exp \left( -\sum_{j=0}^{n-1} r(t_j,\widetilde X_j) \, (t_{j+1} - t_j) \right)
    = \widetilde R_{n-1} \, \exp \big( -r(t_{n-1},\widetilde X_{n-1}) \, (t_{n} - t_{n-1}) \big),
	\qquad
	\widetilde R_0 := 1.
\]
Here, $\widetilde R_n$ is a discrete approximation of the term $\exp \big( - \int_0^{t_n} r(\nu, X_\nu)\, \d \nu \big)$.
In the case of a space-independent or even constant potential function $r(t,x)$, this can be simplified, e.g., $\widetilde R_n = e^{- r \, t_n}$ in the case of a constant potential $r(t,x) = r$.
The discrete approximation~\eqref{eq:YProcessFull} can then be used to generate training samples $\{(x^i,y^i)\}_{i=1}^{n_\text{data}}$ and train a neural network $u_\theta: \DD \to \R$ which approximates the solution of the PDE~\eqref{eq:LinearParabolicPDE} in the domain of interest $\DD$ at time $t=0$.

An alternative formulation of \eqref{eq:FeynmanKacGeneral} can be obtained by means of the concept of \emph{killed} stochastic processes, see  \cite[Sec.\ 8.2]{Oksendal2003} or \cite[Ch.\ 15]{steele2012stochastic}.
Such a process $(\widehat X_t)_{t\in[0,T]}$ behaves exactly like the process $\{X_t\}_{t\in[0,T]}$, but becomes undefined or ``killed'' at a certain random (killing) time $\zeta$, after which the process $\widehat X_t$ is assigned a so-called ``coffin state''.
Here, $\zeta$ is an exponentially distributed random time with ``killing rate'' $r(t,x)$.
Thus, it can be shown, see \cite[Sec.\ 8.2]{Oksendal2003}, that
the solution of the parabolic PDE \eqref{eq:LinearParabolicPDE} admits the representation
\begin{equation} \label{eq:FeynmanKacGeneralKilling}
    u(t,x) 
    = 
    \E \left[
        \int_t^T f(\tau,\widehat X_\tau) \, \d \tau
        + 
        g(\widehat X_T) \given \widehat X_t = x
        \right]
        \qquad
        \forall (t,x) \in [0,T] \times \R^d.
\end{equation}

Finally, we mention that boundary conditions can be incorporated into the PDE-SDE framework by considering certain kinds of stochastic processes.
For example, in the case of a linear parabolic PDE as in~\eqref{eq:LinearParabolicPDE} but posed on a bounded spatial domain $\OO$ in place of $\R^d$, the appropriate concept is that of \emph{stopped processes}, which evolve according to the SDE~\eqref{eq:SDEspec} in $\OO$ and are stopped as soon as they hit the parabolic boundary $(0,T) \times \partial \OO \cup \{T\} \times \overline{\OO}$ where $\overline{\OO}$ denotes the closure of $\OO$.
For further details, see \cite{pardoux1985discretization,milstein1997weak,darling1997backwards,bouchard2004discrete,yang2018first} and the references therein.

\subsection{Summary and Extensions}

The approach discussed in this section can be used to solve backward Kolmogorov equations in high-dimensions.
It is based on the Feynman-Kac connection between SDEs and PDEs and can be implemented efficiently using TensorFlow and other scientific machine learning software environments without deeper knowledge since it reduces, in essence, to a regression problem where the data is sampled either directly or via SDE time-stepping methods such as the Euler-Maruyama scheme.

In \cite{berner2020numerically}, a similar technique is proposed for the solution of \emph{parametric} linear Kolmogorov PDEs.
Again, this methodology generates training data by sampling; the employed neural networks, however, are based on a multilevel architecture with residual connections.

\section{Semilinear PDEs in high Dimensions} \label{sec:DeepBSDESolver}

In this section we extend the methodology of \cref{sec:FeynmanKacSolver} to solving \emph{semilinear} PDEs obtained by allowing the lower-order terms in \eqref{eq:KolmogorovPDE} and \eqref{eq:LinearParabolicPDE} to depend nonlinearly on the solution and its gradient.
This results in the final value problem
\begin{equation} \label{eq:SemilinPDE}
\begin{aligned}
   \partial_t u(t,x) 
   + \frac{1}{2} \sigma \sigma^T(t,x) : \nabla^2 u(t,x) 
   + \mu(t,x) \cdot \nabla u(t,x) 
   + f(t,x,u(t,x),\sigma^T(t,x) \nabla u(t,x)) 
   &= 0, 
   	\quad &&(t,x) \in [0,T) \times \R^d,\\
   u(T,x) &= g(x), \quad &&x \in  \R^d,
\end{aligned}
\end{equation}
with drift $\mu$, diffusion $\sigma$ and final data $g$ as before.
The function $f\colon [0,T] \times \R^d \times \R \times \R^d \to \R$ containing lower order terms can depend in a general way on the independent variables $t,x$ as well as on the solution $u(t,x)$  and its transformed gradient $(\sigma^T\nabla) u(t,x)$.
The non-divergence form of the leading-order term as well as the specific dependence on $\sigma^T\nabla u$ again result from the connection between PDEs and stochastic processes.
As we will see in~\cref{sec:DeepBSDEBackground}, the presence of these dependencies requires extending the numerical solution method to include additional approximating stochastic processes for $\nabla u$.

Problems of the form~\eqref{eq:SemilinPDE} arise 
in physics in the form of, e.g., the Allen-Cahn, Burgers or reaction-diffusion equations; 
in finance, e.g., for pricing derivatives with default risk \cite{el1997backward,brigo2013counterparty,crepey2015bilateral};
and stochastic control problems, see \cite{Pham2009}.
The method discussed below is an extension to that presented in~\cref{sec:FeynmanKacSolver} in that it is also based on the PDE-SDE connection, but in this case it is the correspondence of nonlinear PDEs with \emph{backward stochastic differential equations (BSDEs)} \cite{RevuzYor1999,karatzas2014brownian}.
In the linear case discussed in \cref{sec:FeynmanKacSolver} the approximation of the solution $u$ at time $t=0$ is based on a neural network approximation of the mapping $u(0,\cdot):\DD \to \R$, the Feynman-Kac representation $u(0,x) = \E[g(X_T) \given X_0 = x]$ for $x \in \DD$ and generating a large number of 
sample paths of the stochastic process $\{X_t\}_{t \in [0,T]}$ determined by \eqref{eq:SDE} to approximate the conditional expection and train the model.
Using the theory of BSDEs, it is possible to treat nonlinearities of the type contained in~\eqref{eq:SemilinPDE}.

The specific method presented here was proposed in~\cite{EHanJentzen2017,HanJentzenE2018} and is based on earlier work~\cite{han2016deep}.
Again, the focus lies on solving high-dimensional problems and overcoming one source of the curse of dimensionality \cite{Bellman1957}: a high-dimensional state space (large $d$).
In recent years, a number of approaches have been proposed for mitigating or overcoming the curse of dimensionality in solving high-dimensional PDEs.
In the meantime, a number of theoretical results indicate this may indeed be possible; an (incomplete) list is given in \cref{sec:FurtherMethods}. 
In \cite{hutzenthaler2020proof} it is proven that deep ReLU networks, i.e., neural networks with multiple hidden layers and the rectified linear unit activation function, are in theory able to overcome the curse of dimensionality for certain kinds of the semilinear parabolic equations with nonlinearities which do not involve the gradient.
This is similar to the linear case~\cite{jentzen2018proof}.
In particular, it can be shown that the number of parameters in the neural network grows at most polynomially in both the dimension of the PDE ($d + 1$) and the reciprocal of the desired approximation accuracy.
Note however, that training a neural network in general is a NP-hard problem, \cite[Sec.\ 20.5]{shalev2014understanding}.
The proof relies on full history recursive multilevel Picard approximations, see also~\cite{eEtAl2019,beck2020overcoming}

The approach discussed below can be used to construct an approximate solution of the semilinear problem~\eqref{eq:SemilinPDE} at a fixed point in time over a bounded domain of interest $\mathcal{D} \subset \R^d$ by sampling the initial point $x$ uniformly on $\DD$ as in~\cref{sec:FeynmanKacMethodology}.
For simplicity, however, we consider the problem of determining the solution at a specific point in space and time, i.e., to determine $u(0,x)$ for fixed $x \in \R^d$.

\subsection{Theoretical Background} \label{sec:DeepBSDEBackground}

As in \cref{sec:FeynmanKacSolver}, we consider a time-evolution $\{X_t\}_{t \in [0,T]}$ in state space $\mathbb R^d$ driven by the  \emph{forward} SDE
\begin{equation}  \label{eq:ForwardBSDE}
   X_t = x + \int_0^t \mu(s,X_s) \,\d s + \int_0^t \sigma(s,X_s) \,\d W_s
\end{equation}
starting at $x \in \mathbb R^d$, with underlying probability space $(\Omega, \FF, \P; \F)$ with filtration $\F = \{\FF_t\}_{t \in [0,T]}$ induced by a $d$-dimensional Brownian motion $\{W_t\}_{t \in [0,T]}$.
In \cref{sec:FeynmanKac} we concluded from \Ito's formula in~\eqref{eq:uDyn1}--\eqref{eq:uDyn2} that, given a sufficiently smooth function $v \colon [0,T] \times \mathbb R^d \to \mathbb R$, the dynamics of the \emph{value process} $Y_t := v(t,X_t)$ is governed by the SDE (now written in differential notation) 
\begin{align} \label{valueProcessBSDE}
    \d Y_t &= \left( \partial_t v
    + \frac{1}{2} \sigma \sigma^T : \nabla^2 v
    + \nabla v \cdot \mu
    \right) (t, X_t) \, \d t
    + \big(\sigma^T \nabla v\big) (t,X_t) \cdot  \d W_t.
\end{align}
As in \cref{sec:FeynmanKac} we now assume a sufficiently smooth solution $u$ of \eqref{eq:SemilinPDE} to exist, set $v=u$ in \eqref{valueProcessBSDE}, and introduce a third stochastic process $Z_t :=\big(\sigma^T \nabla u \big)(t, X_t)$ to obtain 
\[
    \d Y_t = -f(t, X_t, Y_t, Z_t) \, \d t + Z_t \cdot \d W_t,
    \qquad 
    Y_T = g(X_T).
\]
This SDE with final condition $Y_T = g(X_T)$ inherited from \eqref{eq:SemilinPDE} is known as the \emph{BSDE associated with}~\eqref{eq:SemilinPDE} and reads, in integral notation, as
\begin{equation} \label{eq:BackwardBSDE}
    Y_t = g(X_T) + \int_t^T f(s, X_s, Y_s, Z_s) \,\d s - \int_t^T Z_s \cdot \,\d W_s.
\end{equation}
Under suitable regularity assumptions on the functions $\mu, \sigma, f$ and $g$, the SDEs \eqref{eq:ForwardBSDE} and \eqref{eq:BackwardBSDE} possess a unique solution $(X_t, Y_t, Z_t)$ and the link to the nonlinear PDE is given by a generalization of the Feynman-Kac formula which states that for all $t \in [0,T]$ there holds $\P$-a.s. that
\begin{equation} \label{eq:ConnectionYandu}
    Y_t = u(t,X_t) \quad \text{and} \quad Z_t = \big(\sigma^T \, \nabla u) (t, X_t).
\end{equation}
In view of the analogy to \eqref{eq:ValueU} and~\eqref{eq:FeynmanKacGeneral}, the identities \eqref{eq:ConnectionYandu} are sometimes referred to as the \emph{nonlinear Feynman-Kac representation} \cite[Sec.\ 6.3]{Pham2009}.
The system consisting of \eqref{eq:ForwardBSDE} and \eqref{eq:BackwardBSDE} is called a \emph{forward-backward stochastic differential equation (FBSDE)} \cite{pardoux1990adapted,pardoux1992backward,pardoux1999forward}.
We note that the forward SDE~\eqref{eq:ForwardBSDE} does not depend on $Y_t$ and $Z_t$, and can thus be solved independently.
As a result, the desired solution value $u(0,x)$ can now be found by solving the FBSDE and evaluating $Y_0$ in ~\eqref{eq:ConnectionYandu}.
We refer to \cite[Ch.\ 7]{yong1999stochastic} for a general account on the solvability of FBSDEs.

The difference to the procedure described in~\cref{sec:FeynmanKacGeneral} is that the solution of the value process $\{Y_s\}_{s \in [0,T]}$ is now more involved due to the nonlinear term $f$ and its dependence on $u(t,x)$ and $(\sigma^T \nabla) u(t,x)$.
% In particular, $Y_s$ appears now on both sides of the equation, and a solution by simple integration is no longer feasible.\\
% \oge{TODO: Diesen Satz verstehe ich immer noch nicht, auf beiden Seiten welcher Gleichung? Gilt das in derselben Weise nicht auch für $X_t$ in \eqref{eq:ForwardBSDE}?}

\subsection{Deep BSDE Solver} \label{sec:DeepBSDEMethodology}

The algorithm termed \emph{deep BSDE solver} in \cite{HanJentzenE2018} constructs an approximation to a solution value $u(0,x)$ of the PDE~\eqref{eq:SemilinPDE} by way of solving the associated FBSDE \eqref{eq:ForwardBSDE}, \eqref{eq:BackwardBSDE}, yielding $u(0,x) = Y_0$ as summarized in \cref{sec:DeepBSDEBackground}.
We now proceed to show how this is achieved using neural networks.

Starting with a discretization of the time domain $[0,T]$ into $N$ equidistant intervals with steps $0=t_0 < t_1 < \ldots < t_N = T$ and step size $\Delta t = T/N$, we generate approximate sample paths of the continuous time process $\{X_t\}_{t \in [0,T]}$ using the Euler-Maruyama scheme for the forward SDE \eqref{eq:ForwardBSDE} which yields the discrete time process
\begin{equation} \label{eq:EMForward}
    \widetilde X_{n+1} 
    = 
    \widetilde X_n + \mu(t_n, \widetilde X_n) \, (t_{n+1} - t_n) 
                   + \sigma(t_n, \widetilde X_n) \, (W_{t_{n+1}} - W_{t_n})
    \quad \text{with} \quad \widetilde X_0 = x.
\end{equation}
In the same way, we construct sample paths for the backward SDE \eqref{eq:BackwardBSDE} as
\begin{equation} \label{eq:EMBackward}
    \widetilde Y_{n+1} 
    = 
    \widetilde Y_n - f(t_n, \widetilde X_n, \widetilde Y_n, \widetilde Z_n) \, (t_{n+1} - t_n)
    			   + \widetilde Z_{n} \cdot (W_{t_{n+1}} - W_{t_n})
    \quad \text{with} \quad \widetilde Y_{N} = g(\widetilde X_{N}).
\end{equation}
Note that the increments of the Brownian motion $(W_{t_{n+1}} - W_{t_n})\sim \pdist{N}(0,(t_{n+1} - t_n) I_{d\times d})$ are the same in~\eqref{eq:EMForward} and~\eqref{eq:EMBackward}.

The algorithm can be summarized by the following steps:
\begin{enumerate}[(1)]
\item
    Simulate paths of the discrete state space process $\{\widetilde X_{n}\}_{n=0}^N$ and the corresponding increments of the Brownian motion $\{W_{t_{n+1}} - W_{t_n}\}_{n=0}^{N-1}$ according to the time-stepping scheme \eqref{eq:EMForward}.
\item
    Simulate paths of the discrete value process $\{\widetilde Y_{n}\}_{n=0}^N$ according to the time-stepping scheme ~\eqref{eq:EMBackward}.
Closer inspection reveals that \eqref{eq:EMBackward} contains unknown quantities necessary to carry out the time-stepping: $\widetilde Y_{0}$, which is an approximation of $u(0, x)$ as well as $\widetilde Z_{n}$ for $n = 0,\ldots,N-1$, which are approximations of $(\sigma^T\nabla u)(t_n, \widetilde X_{n})$.
These quantities are obtained by training a neural network.

The quantities $\widetilde Y_{0} \approx u(0, x)$ and $\widetilde Z_{0} \approx (\sigma^T\nabla u)(0, x)$ are treated as individual parameters---both needed only in the point $(0,x)$---and are learned in the course of training.
The remaining quantities $\widetilde Z_{n}, n=1,\dots,N-1$ are approximated by neural networks which realize the mapping $x \mapsto (\sigma^T\nabla u)(t_n, x)$ for $n=1,\ldots,N-1$.
All neural network parameters to be learned are collected in
\[
    \theta = \left( \theta_{u_0}, \theta_{\nabla u_0}, \theta_{\nabla u_1}, \ldots, \theta_{\nabla u_{N-1}} \right),
\]
        where $\theta_{u_0} \in \R$, $\theta_{\nabla u_0} \in \R^d$  and $\theta_{\nabla u_n} \in \R^{\rho_n}$ and $\rho_n$ is the number of unknown parameters in the neural network realizing the mapping $x \mapsto (\sigma^T\nabla u)(t_n,x)$ for $n=1\ldots,N-1$.
\item
Since $\widetilde Y_{N}$ should approximate $u(T, \widetilde X_{N}) = g(\widetilde X_{N})$ according to~\eqref{eq:EMBackward} the network is trained to minimize the mean squared error (MSE) between $\widetilde Y_{N}$ and $g(\widetilde X_{N})$.
For a batch of $m$ simulated pairs $(\widetilde X_{N}, \widetilde Y_{N})$, this results in the loss function 
\[
   \phi_\theta(\widetilde X_N, \widetilde Y_N) 
   := 
   \frac{1}{m}\sum_{i=1}^m \left[\widetilde Y^i_{N} - g(\widetilde X^i_{N})\right]^2,
            % \phi(X,Y) := \frac{1}{m}\sum_{i=1}^m \left(Y^i_T - g(X^i_T)\right)^2.
\]
where $Y^i_N$ is the output of the neural network.
        Automatic differentiation of $\phi_\theta$ with respect to the unknowns $\theta$ is then employed to obtain the gradient $\nabla_\theta \phi_\theta$, which is then used by an optimization routine, e.g., some variant of the stochastic gradient descent method.%
        Note that the same considerations with regard to overfitting as noted at the end of \Cref{sec:FeynmanKacTraining} in connection with the Feynman-Kac solver apply here.
        
\end{enumerate}

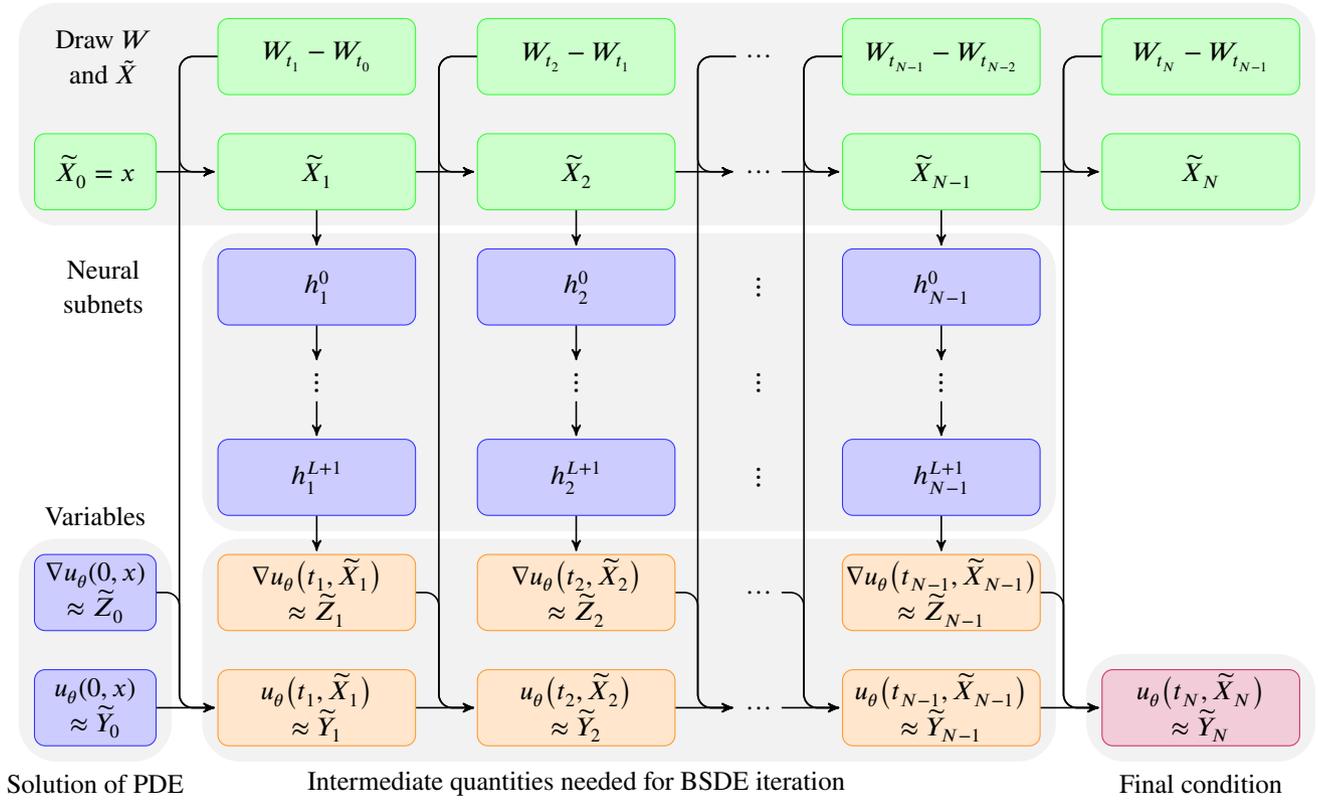
\begin{figure}[htbp]
\centering
\tikzset{%
    cell/.style={rounded corners,
        minimum width=2.6cm,
        minimum height=1.0cm,
        inner sep=0.04cm,
        font=\normalsize},
    Terminalstyle/.style={cell,
        draw=purple!80,
        fill=purple!20},
    Xstyle/.style={cell,
        draw=green!80,
        fill=green!20},
    Ystyle/.style={cell,
        draw=orange!80,
        fill=orange!20},
    Zstyle/.style={cell,
        draw=orange!80,
        fill=orange!20},
    hstyle/.style={cell,
        draw=blue!80,
        fill=blue!20},
    gradstyle/.style={cell,
        draw=orange!80,
        fill=orange!20},
    vstyle/.style={cell,
        draw=red!80,
        fill=red!20},
    pstyle/.style={->,
        shorten >=1pt,
        >=stealth',
        semithick,
        rounded corners=5pt},
    background/.style={rectangle,
        fill=gray!10,
        inner sep=0.2cm,
        rounded corners=5mm},
    varbackground/.style={rectangle,
        fill=gray!10,
        inner sep=0.2cm,
        rounded corners=5mm}
    }

    \begin{tikzpicture}[>=latex]
        % ,text height=1.0ex,text depth=0.25ex]
    % "text height" and "text depth" are required to vertically
    % align the labels with and without indices.
  
  % The various elements are conveniently placed using a matrix:
  \matrix[row sep=0.5cm,column sep=0.8cm] {
    % First line: dW
      &
      \node (dw_0) {}; &
        \node (dw_1) [Xstyle]{$W_{t_1} - W_{t_0}$}; &
        \node (dw_2) [Xstyle]{$W_{t_2} - W_{t_1}$}; &
        \node (dw_p) {$\ldots$}; &
        \node (dw_Nm1) [Xstyle]{$W_{t_{N-1}} - W_{t_{N-2}}$}; &
        \node (dw_N) [Xstyle]{$W_{t_N} - W_{t_{N-1}}$}; &
        \\
        &
        \node (X_0) [Xstyle,minimum width=1.6cm]{$\widetilde X_0 = x$}; &
        \node (X_1) [Xstyle]{$\widetilde X_{1}$}; &
        \node (X_2) [Xstyle]{$\widetilde X_{2}$}; &
        \node (X_p) {$\ldots$}; &
        \node (X_Nm1) [Xstyle]{$\widetilde X_{{N-1}}$}; &
        \node (X_N) [Xstyle]{$\widetilde X_{N}$}; &
        \\
        & & 
        \node (h1_1) [hstyle]{$h_1^0$}; &
        \node (h1_2) [hstyle]{$h_2^0$}; &
        \node (h1_p) {$\vdots$}; &
        \node (h1_Nm1) [hstyle]{$h_{N-1}^0$}; &
        % \node (h1_N) [hstyle]{$h_N^1$}; &
        \\
        & & 
        \node (hp_1) {$\vdots$}; &
        \node (hp_2) {$\vdots$}; &
        \node (hp_p) {$\vdots$}; &
        \node (hp_Nm1) {$\vdots$}; &
        % \node (hp_N)  {$\vdots$}; &
        \\
        & & 
        \node (hL_1) [hstyle]{$h_1^{L+1}$}; &
        \node (hL_2) [hstyle]{$h_2^{L+1}$}; &
        \node (hL_p) {$\vdots$}; &
        \node (hL_Nm1) [hstyle]{$h_{N-1}^{L+1}$}; &
        % \node (h1_1) [hstyle]{$h_1^1$}; &
        % \node (h1_2) [hstyle]{$h_2^1$}; &
        % \node (h1_p) {$\vdots$}; &
        % \node (h1_Nm1) [hstyle]{$h_{N-1}^1$}; &
        % % \node (h1_N) [hstyle]{$h_N^1$}; &
        % \\
        % & & 
        % \node (hp_1) {$\vdots$}; &
        % \node (hp_2) {$\vdots$}; &
        % \node (hp_p) {$\vdots$}; &
        % \node (hp_Nm1) {$\vdots$}; &
        % % \node (hp_N)  {$\vdots$}; &
        % \\
        % & & 
        % \node (hL_1) [hstyle]{$h_1^L$}; &
        % \node (hL_2) [hstyle]{$h_2^L$}; &
        % \node (hL_p) {$\vdots$}; &
        % \node (hL_Nm1) [hstyle]{$h_{N-1}^L$}; &
        % % \node (hL_N) [hstyle]{$h_N^L$}; &
        \\
        &
        \node (grad_0) [hstyle,align=center,minimum width=1.6cm]{$\nabla u_\theta(0,x)$\\$\approx \widetilde Z_0$}; &
        \node (grad_1) [Zstyle,align=center]{$\nabla u_\theta\big(t_1,\widetilde X_{1}\big)$\\$\approx \widetilde Z_1$}; &
        \node (grad_2) [Zstyle,align=center]{$\nabla u_\theta\big(t_2,\widetilde X_{2}\big)$\\$\approx \widetilde Z_2$}; &
        \node (grad_p) {$\ldots$}; &
        \node (grad_Nm1) [Zstyle,align=center]{$\nabla u_\theta\big(t_{N-1},\widetilde X_{{N-1}}\big)$\\$\approx \widetilde Z_{N-1}$}; &
        \node (grad_N) [Zstyle,align=center,draw=none,fill=none]{\phantom{$Z_{N-1}\approx$}\\\phantom{$\nabla u_\theta\big(t_{N},\widetilde X_{{N}}\big)$}}; &
        % Only nabla u
        % \node (grad_0) [hstyle,minimum width=1.6cm]{$\nabla u_\theta\big(t_0,\xi\big)$}; &
        % \node (grad_1) [Xstyle]{$\nabla u_\theta\big(t_1,\widetilde X_{1}\big)$}; &
        % \node (grad_2) [Xstyle]{$\nabla u_\theta\big(t_2,\widetilde X_{2}\big)$}; &
        % \node (grad_p) {$\ldots$}; &
        % \node (grad_Nm1) [Xstyle]{$\nabla u_\theta\big(t_{N-1},\widetilde X_{{N-1}}\big)$}; &
        % \node (grad_N) [Xstyle,draw=none,fill=none]{\phantom{$\nabla u_\theta\big(t_{N},\widetilde X_{{N}}\big)$}}; &
        \\
        &
        \node (u_0) [hstyle,align=center,minimum width=1.6cm]{$u_\theta(0,x)$\\$\approx\widetilde Y_0$}; &
        \node (u_1) [Ystyle,align=center]{$u_\theta\big(t_1,\widetilde X_{1}\big)$\\$\approx\widetilde Y_1$}; &
        \node (u_2) [Ystyle,align=center]{$u_\theta\big(t_2,\widetilde X_{2}\big)$\\$\approx\widetilde Y_2$}; &
        \node (u_p) {$\ldots$}; &
        \node (u_Nm1) [Ystyle,align=center]{$u_\theta\big(t_{N-1},\widetilde X_{{N-1}}\big)$\\$\approx\widetilde Y_{N-1}$}; &
        \node (u_N) [Terminalstyle,align=center]{$u_\theta\big(t_{N},\widetilde X_{{N}}\big)$\\$\approx\widetilde Y_N$}; &
        % \node (u_N) [Terminalstyle,align=center]{$u_\theta\big(t_{N},\widetilde X_{{N}}\big)$\\$\approx\widetilde Y_N=g\big(\widetilde X_N\big)$}; &
        % \node (u_0) [hstyle,minimum width=1.6cm]{$u_\theta\big(t_0,\xi\big)$}; &
        % \node (u_1) [Xstyle]{$u_\theta\big(t_1,\widetilde X_{1}\big)$}; &
        % \node (u_2) [Xstyle]{$u_\theta\big(t_2,\widetilde X_{2}\big)$}; &
        % \node (u_p) {$\ldots$}; &
        % \node (u_Nm1) [Xstyle]{$u_\theta\big(t_{N-1},\widetilde X_{{N-1}}\big)$}; &
        % \node (u_N) [Terminalstyle]{$u_\theta\big(t_{N},\widetilde X_{{N}}\big)$}; &
        \\
    };

    % Connections X -> X
    \draw[pstyle] (X_0) -- (X_1);
    \draw[pstyle] (X_1) -- (X_2);
    \draw[pstyle] (X_2) -- (X_p);
    \draw[pstyle] (X_p) -- (X_Nm1);
    \draw[pstyle] (X_Nm1) -- (X_N);
    % Connections dW -> X
    \draw[pstyle] (dw_1.west) -- ++(-0.5,0.) |- (X_1.west);
    \draw[pstyle] (dw_2.west) -- ++(-0.5,0.) |- (X_2.west);
    \draw[pstyle] (dw_p.west) -- ++(-0.5,0.) |- (X_p.west);
    \draw[pstyle] (dw_Nm1.west) -- ++(-0.5,0.) |- (X_Nm1.west);
    \draw[pstyle] (dw_N.west) -- ++(-0.5,0.) |- (X_N.west);
    % Connections X -> h^1
    \draw[pstyle] (X_1) -- (h1_1);
    \draw[pstyle] (X_2) -- (h1_2);
    \draw[pstyle] (X_Nm1) -- (h1_Nm1);
    % Connections h^1 -> h^p
    \draw[pstyle] (h1_1)   -- (hp_1)  ;
    \draw[pstyle] (h1_2)   -- (hp_2)  ;
    \draw[pstyle] (h1_Nm1) -- (hp_Nm1);
    % Connections h^p -> h^L
    \draw[pstyle] (hp_1)   -- (hL_1)  ;
    \draw[pstyle] (hp_2)   -- (hL_2)  ;
    \draw[pstyle] (hp_Nm1) -- (hL_Nm1);
    % Connections h^L -> grad
    \draw[pstyle] (hL_1)   -- (grad_1)  ;
    \draw[pstyle] (hL_2)   -- (grad_2)  ;
    \draw[pstyle] (hL_Nm1) -- (grad_Nm1);
    % Connections u -> u
    \draw[pstyle] (u_0) -- (u_1);
    \draw[pstyle] (u_1) -- (u_2);
    \draw[pstyle] (u_2) -- (u_p);
    \draw[pstyle] (u_p) -- (u_Nm1);
    \draw[pstyle] (u_Nm1) -- (u_N);
    %Connections W -> u
    \draw[pstyle] (dw_1.west) -- ++(-0.5,0.) |- (u_1.west);
    \draw[pstyle] (dw_2.west) -- ++(-0.5,0.) |- (u_2.west);
    \draw[pstyle] (dw_p.west) -- ++(-0.5,0.) |- (u_p.west);
    \draw[pstyle] (dw_Nm1.west) -- ++(-0.5,0.) |- (u_Nm1.west);
    \draw[pstyle] (dw_N.west) -- ++(-0.5,0.) |- (u_N.west);
    % Connections grad_{k-1} -> u_k
    \draw[pstyle] (grad_0.east) -- ($(grad_1.west) - (0.5,0.)$) |- (u_1.west);
    \draw[pstyle] (grad_1.east) -- ($(grad_2.west) - (0.5,0.)$) |- (u_2.west);
    \draw[pstyle] (grad_2.east) -- ($(grad_p.west) - (0.5,0.)$) |- (u_p.west);
    \draw[pstyle] (grad_p.east) -- ($(grad_Nm1.west) - (0.5,0.)$) |- (u_Nm1.west);
    \draw[pstyle] (grad_Nm1.east) -- ($(grad_N.west) - (0.5,0.)$) |- (u_N.west);

    \begin{pgfonlayer}{background}
        \node [background,
                    fit=(X_0) (dw_N)] {};
        \node [varbackground,
                    fit=(h1_1) (hL_Nm1)] {};
        \node [varbackground,
                    fit=(u_0) (grad_0)] {};
        \node [background,
                    fit=(grad_1) (u_Nm1)] {};
        \node [background,
                    fit=(u_N)] {};
    \end{pgfonlayer}

    % Custom labels
    \node (lab_var) at ($(grad_0) + (0.0,1.0)$) {Variables};
    \node[align=center] (lab_subnets) at ($(h1_1.west) - (1.5,0.0)$) {Neural\\subnets};
    \node[align=center] (lab_X) at ($(dw_1.west) - (1.5,0.0)$) {Draw $W$\\and $\tilde X$};
    \node (lab_qoi) at ($(u_0) - (0.0,1.0)$) {Solution of PDE};
    \node[align=center] (lab_intermed) at ($(u_2) - (0.0,1.0)$) {Intermediate quantities needed for BSDE iteration};
    \node[align=center] (lab_terminal) at ($(u_N) - (0.0,1.0)$) {Final condition};

\end{tikzpicture}

\caption{Illustration of the complete deep BSDE solver model adapted from~\cite{EHanJentzen2017,HanJentzenE2018} in the case $\sigma = I_{d \times d}$.
The two upper rows express the evolution of the forward process $\{X\}_{n=0}$ starting at $\widetilde X_0 = x$ (\emph{green}).
The unknown parameters for $u_\theta (0,x)$ and $\nabla u_\theta(0,x)$ (\emph{blue, left}) as well as the parameters in the neural network approximating $\widetilde Z_n$, $n=1,\ldots,N-1$ (\emph{blue, center}) are learned by training.
The intermediate values $\widetilde Y_n$ and $\widetilde Z_n$, $n=1,\ldots,N-1$ (\emph{orange}) are needed to establish the link between the desired PDE solution value $u(0,x) \approx \widetilde Y_0$ with the given final value $\widetilde Y_N = g\big(\widetilde X_N\big)$ (\emph{red}).
} \label{fig:DeepBSDEModel}
\end{figure}

The complete network structure is illustrated in \cref{fig:DeepBSDEModel}.
The architecture of the sub-networks realizing the mapping $x \mapsto (\sigma^T\nabla u)(t_n,x)$ used in the numerical experiments described below are taken to be the same as in ~\cite{HanJentzenE2018}, where they are given by
\begin{equation} \label{eq:DeepBSDENetwork}
	\texttt{Input \lto BN \lto (Dense \lto BN \lto ReLU) \lto (Dense \lto BN \lto ReLU) \lto Dense \lto BN
	\lto Output}
\end{equation}
Here, \texttt{BN} stands for \emph{batch normalization},
\texttt{Dense} denotes a fully connected layer without bias term and activation, and
\texttt{ReLU} denotes the application of the componentwise rectified linear unit activation function $\relu(x) = \max \{0,x\}$.
In terms of the layers in \cref{fig:DeepBSDEModel} this means the following: 
first, the inputs $\widetilde X_n \in \mathbb R^d$ are scaled and shifted componentwise by batch normalization, resulting in $h^1_n := \texttt{BN}^1_n(\widetilde X_n)$;
second, the outputs from the first layer are processed by the subsequent block
$h^2_n := \relu\big( \bn^2_n(W^2_n h^1_n)\big)$
followed by block
$h^3_n := \relu\big( \bn^3_n(W^3_n h^2_n)\big)$;
finally, the output is multiplied by another matrix $W^4_n$ and batch normalized once more, giving
$h^4_n := \bn^4_n(W^4_n h^3_n) \approx \widetilde Z_n$.

To implement the model in TensorFlow~\cite{tensorflow2015-whitepaper} all that is needed is to provide a routine that realizes the interaction between the known and unknown quantities and respects the time-stepping scheme  \eqref{eq:EMBackward}.
In the following, we discuss two examples.
An implementation of the methodology for both examples is given in the accompanying Jupyter notebook \texttt{DeepBSDE\_Solver.ipynb}.

\subsection{Example: Linear-Quadratic Gaussian Control} \label{sec:LinearQuadraticGaussianControl}

We consider the linear-quadratic Gaussian control problem as discussed in \cite[Sec.\ 4.3]{EHanJentzen2017}, 
\cite{HanJentzenE2018} and 
\cite{bachouch2020deep}.
The goal is to control a stochastic process $\{X_t\}_{t \in [0,T]}$ governed by the SDE
\begin{align*}
	X_t = x +  2 \int_0^t  m_s \, \d s + \sqrt{2} \int_0^t  \d W_s
	% X_t = x +  2 \sqrt{\lambda} \int_0^t  m_s \, \d s + \sqrt{2} \int_0^t  \d W_s
\end{align*}
with a control $m_t \in \R^d$ entering as the drift term.
The solution of the control problem is characterized by the value function, i.e., the function $u\colon [0,T] \times \R^d \to \R$ that gives the minimal expected sum of accumulated running cost and final cost over all admissible control processes\footnote{%
In this setting, an $\mathbb R^d$-valued control process $\{m_s\}_{s \in [t,T]}$ is admissible if its value at time $s$ is based only on the information available up to time $s$. 
To be precise, the process $m_s$ has to be progressively measurable with respect to the underlying filtration $\F$; see \cite{Pham2009,yong1999stochastic} for further details.
} from time $t$ onward starting at $x$:
\begin{align}
    \label{eq:HJBValue}
    u(t,x) 
    = 
    \min_{\{m_s\}_{s \in [t,T]}} 
    \mathbb{E} \left[ \int_t^T \|m_s\|^2 \, \d s + g(X_T) \given X_t = x\right]
    % u(t,x) = \min_m \mathbb{E} \left[ \int_t^T \|m_t\|^2 \mathrm{d}t + g(X_T) \given X_t = x\right]
\end{align}
The function $g:\mathbb{R}^d \to \mathbb{R}$ is the prescribed final data.
% and $\lambda > 0$ is a fixed parameter describing the intensity of the control.
The Hamilton-Jacobi-Bellman equation associated with the stochastic control problem is given by the nonlinear PDE
\begin{align} \label{eq:HJBLQG}
\begin{aligned}
    \partial_t u(t,x) + \Delta u(t,x) 
    + \min_m \left\{ 2 \, m \cdot \nabla u(t,x)  + \|m\|^2\right\} 
    % + \min_m \left\{ 2 \sqrt{\lambda} \, m \cdot \nabla u(t,x)  + \|m\|^2\right\} 
    &= 0,  	\quad &&(t,x) \in [0,T) \times \R^d,\\
    u(T,x) &= g(x), && x \in \mathbb R^d.
\end{aligned}
\end{align}
Note that this equation is purely deterministic.
As easily verified, the minimum is attained at $m = - \nabla u$.
% As easily verified, the minimum is attained at $m = - \sqrt{\lambda} \nabla u$.
Inserting this optimal control into the HJB equation~\eqref{eq:HJBLQG} yields the semilinear PDE
\begin{align} \label{eq:GaussControlPDE}
\begin{aligned}
    \partial_t u(t,x) + \Delta u(t,x) - \|\nabla u(t,x)\|^2 &= 0, \quad &&(t,x) \in [0,T) \times \R^d,\\
    u(T,x) &= g(x), && x \in \mathbb R^d.
\end{aligned}
\end{align}
The formulation~\eqref{eq:HJBValue} reveals that the PDE solution $u$ is the value function of a stochastic control problem,
the control $m_s$ is the negative gradient of the solution $u$ which plays the role of a policy function in a reinforcement learning approach to solve the stochastic control problem \cite{bengio2009learning,lecun2015deep,goodfellow2016deep,sutton2018reinforcement}.
This connection to stochastic control problems provided the original motivation for the deep BSDE method \cite{han2016deep,EHanJentzen2017}.

We solve this equation in dimension $d=100$ for drift coefficient $\mu \equiv 0$, diffusion coefficient $\sigma \equiv \sqrt{2} I_{d \times d}$, reaction term $f(t,x,y,z) = -1/2 \, \norm{z}^2$ and final time $T=1$ 
% We solve this equation in dimension $d=100$ for drift coefficient $\mu \equiv 0$, diffusion coefficient $\sigma \equiv \sqrt{2} I_{d \times d}$, reaction term $f(t,x,y,z) = - \lambda/2 \, \norm{z}^2$ and final time $T=1$ 
with prescribed data
$
    g(x) = \log \big( 1/2 \, (1 + \norm{x}^2) \big)
$
using the algorithm described in~\cref{sec:DeepBSDEMethodology} to approximate the solution value $u(0,x)$ for $x = 0 \in \R^d$.
We note that the solution to this control problem can be obtained explicitly via a Cole-Hopf transformation, see e.g.\ \cite{chassagneux2016numerical}, and is given by the formula
$ u(t,x) = - \log \big(
    \E \big[
        \exp \big(- g\big(x + \sqrt 2 W_{T-t}\big)\big)
        \big]
    \big).
$
% $ u(t,x) = - \lambda^{-1} \log \big(
%     \E \big[
%         \exp \big(- \lambda \, g\big(x + \sqrt 2 W_{T-t}\big)\big)
%         \big]
%     \big).
% $
This can be used as a reference solution.

The results for 4 different experimental configurations are presented in \cref{table:HJBresults}.
All experiments employ the Adam optimizer \cite{kingma2014adam} with 
constant learning rate $\delta = 0.01$ as used in \cite{EHanJentzen2017}, 
the number of training epochs set to $n_\text{epochs}=\num{2000}$ and 
batch size $n_\text{batch}=64$.   
The setup in the second row labeled \emph{Reference} uses the same configuration as employed in~\cite[Sec.\ 4.3]{EHanJentzen2017}, i.e., 
$N = 20$ discrete time steps,
and the network architecture as shown in~\eqref{eq:DeepBSDENetwork} containing two stacks of layers of the form
\begin{equation} \label{eq:DeepBSDEStack}
    \texttt{Dense \lto BN \lto ReLU}
\end{equation}
with $\num{110}$ neurons in each layer.
The \emph{Simple} configuration contains no such layer stack~\eqref{eq:DeepBSDEStack}, and uses only $N=1$ time step, which explains the fast computation.
In setting \emph{L=3}, we increased the number of hidden layer stacks~\eqref{eq:DeepBSDEStack} to three, the number of time steps to $N=30$ and the number of neurons in each layer to $\num{200}$.
In setting \emph{L=5}, we increased the number of hidden layer stacks~\eqref{eq:DeepBSDEStack} to five, the number of time steps to $N=50$ and the number of neurons in each layer to $\num{300}$.

\pgfplotstableset{
    columns/id/.style={
        column name=Experiment,
        string type},
    multicolumn names, % allows to have multicolumn names
    col sep=comma, % the seperator in our .csv file
    sci zerofill,
    columns/L1relMean/.style={
        column name={Mean relative error},
        sci},
    columns/L1relStd/.style={
        column name={Std.-dev.\ relative error},
        sci},
    columns/yMean/.style={
        column name={Mean $u(0,x)$},
        fixed zerofill,
        precision=4},
    columns/yStd/.style={
        column name={Std.-dev.\ $u(0,x)$},
        sci},
    columns/TimeMean/.style={
        column name={Mean time [s]},
        column type={S},string type},
    every head row/.style={
        before row={\toprule}, % have a rule at top
        after row={\midrule}, % rule under units
    },
    % every nth row={3}{before row=\midrule},
    every last row/.style={after row=\bottomrule}, % rule at bottom
    }
    \begin{table}[htb]
        \begin{center}
            \begin{filecontents*}{./table_HJB.csv}
id,Iter,LossMean,LossStd,yMean,yStd,L1relMean,L1relStd,L1absMean,L1absStd,TimeMean,TimeStd
\emph{Simple ($L=0$)},1999,0.02427262877451947,0.0046814579761194275,4.600002347264344,0.0014809137570376248,0.0021469974001172335,0.00032262880820188794,0.009855034346988667,0.0014809137570376248,3.43613862991333,0.0354719858490011
\emph{Reference ($L=2$)},1999,0.02411662391489849,0.0032789221379314754,4.598923819051447,0.0009712422646300817,0.001912031474326298,0.00021159283099626494,0.008776506134092265,0.0009712422646300817,83.91181168556213,6.994217030955384
\emph{$L=3$},1999,0.019164184077715057,0.003537058855354971,4.599102212217818,0.00119326427401168,0.0019508958405891948,0.0002599620867621519,0.008954899300462138,0.00119326427401168,135.77397108078003,12.1501852064883
\emph{$L=5$},1999,0.02229944263148108,0.0038094401219392784,4.598270226369365,0.001720705464550229,0.0017696411244037757,0.00037486933365034027,0.00812291345201004,0.001720705464550229,363.1356254577637,33.664767627904
\end{filecontents*}
\pgfplotstabletypeset[
                columns={{id},{yMean},{yStd},{L1relMean},{L1relStd},{TimeMean}},
            ]{./table_HJB.csv} % filename/path to file
            % id,Iter,LossMean,LossStd,yMean,yStd,L1relMean,L1relStd,L1absMean,L1absStd,TimeMean,TimeStd
        \end{center}
            \caption{Shown are the mean and standard deviations of $u_\theta(0,x)$ and the relative error $\abs{u_\theta(0,x) - u^*}/u^*$, resp., with $u^* \approx 4.5901$ (determined via Monte-Carlo sampling), as well as the mean computation time over 5 consecutive runs with randomly initialized parameters $\theta$ after $n_\text{epochs}=\num{2000}$ training epochs.}
            \label{table:HJBresults}
    \end{table}

	The results in \cref{table:HJBresults} suggests that for the solution of the linear-quadratic Gaussian control problem~\eqref{eq:GaussControlPDE} all models display similar performance.
    It is surprising that even the \emph{Simple} model taking less than $4$ seconds total computation time provides essentially the same approximation quality as the more complex models.
    This is in line with the findings in~\cite{bachouch2020deep} that it appears difficult to further decrease the relative errors using the proposed methodology.
	On the other hand, a decrease in relative approximation error when increasing the number of hidden layers was observed in another example given in \cite{HanJentzenE2018}.
	The convergence behavior of this method seems to call for further research.

\subsection{Example: Allen-Cahn Equation}

As a second example, we solve the Allen-Cahn equation with a double-well potential \cite[Sec.\ 4.2]{EHanJentzen2017}, \cite{HanJentzenE2018,emmerich2003diffuse}, i.e., the semilinear reaction-diffusion equation 
\begin{align*} 
    u_t(t,x) + \Delta u(t,x) + u(t,x) - u^3(t,x) &= 0\\
    u(T,x) &= \big( 2 + \frac{2}{5} \, \norm{x}^2 \big)^{-1}.
\end{align*}
The results for the approximation of $u(0,x)$ for $x = 0 \in \R^d$ with $d=100$ and $T=1$ are displayed in \cref{table:AllenCahnresults}.
The training was carried out over $n_\text{epochs}=\num{4000}$ epochs with the Adam optimizer \cite{kingma2014adam} with a constant step size $\delta = 5 \cdot 10^{-4}$ for the same set of network configurations as used in~\cref{sec:LinearQuadraticGaussianControl}

\pgfplotstableset{
    columns/id/.style={
        column name=Experiment,
        string type},
    multicolumn names, % allows to have multicolumn names
    col sep=comma, % the seperator in our .csv file
    sci zerofill,
    columns/L1relMean/.style={
        column name={Mean relative error},
        sci},
    columns/L1relStd/.style={
        column name={Std.-dev.\ relative error},
        sci},
    columns/yMean/.style={
        column name={Mean $u(0,x)$},
        fixed zerofill,
        precision=6},
    columns/yStd/.style={
        column name={Std.-dev.\ $u(0,x)$},
        sci},
    columns/TimeMean/.style={
        column name={Mean time [s]},
        column type={S},string type},
    every head row/.style={
        before row={\toprule}, % have a rule at top
        after row={\midrule}, % rule under units
    },
    % every nth row={3}{before row=\midrule},
    every last row/.style={after row=\bottomrule}, % rule at bottom
    }
    \begin{table}[htb]
        \begin{center}
            \begin{filecontents*}{./table_AllenCahn.csv}
id,Iter,LossMean,LossStd,yMean,yStd,L1relMean,L1relStd,L1absMean,L1absStd,TimeMean,TimeStd
\emph{Simple ($L=0$)},3999,2.717115699465194e-05,1.0338676244150131e-05,0.055815910040797466,4.5512167996920636e-05,0.05707946745951793,0.0008619402294784412,0.0030139100407974655,4.551216799692063e-05,6.706873369216919,0.05284224328218094
\emph{Reference ($L=2$)},3999,0.00015007395335703118,3.086499834649657e-05,0.05299016853981553,0.0001207164477716472,0.0036288972097957227,0.0021811819981587973,0.00019161303047163374,0.00011517077186678081,157.53529653549194,17.239311158892
\emph{$L=3$},3999,0.0001553455728472791,2.857219767062586e-05,0.052899560662471826,0.0002559106165712257,0.004015599061874192,0.0032830566451238184,0.00021203166166508107,0.00017335195697582787,308.4500723361969,5.461312095502731
\emph{$L=5$},3999,0.0002367898616078146,4.117659396767196e-05,0.052725709399589094,0.00011844214271046397,0.002271589534916321,0.0013996849068049094,0.00011994447062265157,7.390616244911282e-05,637.1849482536315,56.65034158564851
\end{filecontents*}
\pgfplotstabletypeset[
                columns={{id},{yMean},{yStd},{L1relMean},{L1relStd},{TimeMean}},
            ]{./table_AllenCahn.csv} % filename/path to file
            % id,Iter,LossMean,LossStd,yMean,yStd,L1relMean,L1relStd,L1absMean,L1absStd,TimeMean,TimeStd
        \end{center}
            \caption{Shown are the mean and standard deviations of $u_\theta(0,x)$ and the relative error $\abs{u_\theta(0,x) - u^*}/u^*$, resp., with $u^* \approx 0.052802$ (taken from~\cite{EHanJentzen2017}, calculated by a branching-diffusion method), as well as the mean computation time over 5 consecutive runs with randomly initialized parameters $\theta$ after $n_\text{epochs}=\num{4000}$ training epochs.}
            \label{table:AllenCahnresults}
    \end{table}

For this experiment, the difference between the \emph{Simple} and more complex models is clearly visible.
However, the \emph{Simple} model yields again a rough approximation of the solution within only 7 seconds.
Again, the decrease of the relative error is quite small for deeper and wider neural networks with more time steps, similar to our findings in~\cref{sec:LinearQuadraticGaussianControl}.
We note that the accompanying Jupyter notebook \texttt{DeepBSDE\_Solver.ipynb} contains the Burgers-type PDE from \cite[Sec.\ 4.5]{EHanJentzen2017} as a third example.

\subsection{Summary and Extensions}

We have described the deep BSDE solver presented and developed in~\cite{EHanJentzen2017} and~\cite{HanJentzenE2018} for the solution of semilinear PDEs~\eqref{eq:SemilinPDE}.
Note, however, that the solver can also be used to solve BSDEs directly (without taking care of any PDE).

In \cite{BeckEJentzen2019}, the deep BSDE solver considered in this section is extended to fully nonlinear PDEs of second-order.
Here, neural networks are employed to approximate the second-order derivatives of $u$ at a finite number of time steps, from which approximations of the gradients $\nabla u(t_n,\cdot)$ and the function values $u(t_n, \cdot)$ can be derived, similar to~\eqref{eq:EMBackward}.
The method relies on the connection between fully nonlinear second-order PDEs and second-order BSDEs \cite{CheriditoSonerTouziVictoir2007}.

The technique described in~\cite{beck2019deep} is closely related and applies operator splitting techniques to derive a learning approach for the solution of parabolic PDEs in up to $\num{10000}$ spatial dimensions.
In contrast to the deep BSDE method, however, the PDE solution at some discrete time snapshots is approximated by neural networks directly.

Another extension of the deep BSDE solver is considered in~\cite{chanwainam2018machine} where the authors employ a number of adaptations to the proposed methodology in order to improve the convergence properties of the algorithm, e.g., by substituting the activation functions, removing some of the batch normalization layers and using only one instead of $N-2$ neural networks to approximate the scaled gradients of the solution $(\sigma^T \nabla u)(t_n,x)$ for $n=1,\ldots,N-1$.
Furthermore, residual connections are added and more elaborative \emph{long short-term memory (LSTM)} neural networks are employed.
Similarly, the authors in \cite{fujii2019asymptotic} consider the use of asymptotic expansion as prior knowledge in order to improve the accuracy and speed of convergence of the deep BSDE solver.

In~\cite{henry2017deep}, an extension based on a primal-dual solution method for BSDEs using neural networks and a dual formulation of stochastic control problems is discussed, see also \cite{henry2016dual}.
An approach that uses the associated FBSDE to train a neural network to learn the solution of a semilinear PDE is discussed in \cite{RaissiForwardBackward}.

\section{Extensions and Related Work} \label{sec:FurtherMethods}

Beyond the three approaches discussed in detail in \cref{sec:PINNs,sec:DeepBSDESolver,sec:FeynmanKacSolver}, the rapidly developing discipline of scientific machine learning has brought forth a number of promising approaches for solving PDEs beyond the capabilities of conventional numerical methods.
In this final section, we want to give a brief and necessarily incomplete overview over some recent developments.

Before we provide more references concerning neural network-based solution approaches for differential equations we list some results concerning general approximation properties of neural networks.
Early work from the 1990s is now considered foundational, e.g., 
\cite{%
cybenko1989approximation,%
hornik1990universal,%
hornik1991approximation,%
mhaskar1996neural,%
pinkus1999approximation}.
Beginning around 2016, the spectacular successes of machine learning systems in computer vision, natural language processing and other areas prompted renewed efforts to establish a mathematically rigorous foundation for, in particular, deep feedforward neural networks
\cite{%
mhaskar2016deep,%
yarotsky2017error,%
yarotsky2018universal,%
perekrestenko2018universal,%
petersen2018optimal,%
elbrachter2018dnn,%
montanelli2019new,%
bolcskei2019optimal,%
petersen2020topological,%
opschoor2020deep,%
beneventano2020highdimensional,%
gonon2020deep,%
MOPS20_2877,%
GS21_2893,%
laakmann2021efficient,%
HOS21_2887}.
We draw particular attention to a number of publications that rigorously establish that certain neural network architectures are theoretically able to overcome the curse of dimensionality for various linear and nonlinear PDEs, cf.\ \cite{%
jentzen2018proof,%
grohs2018proof,%
hutzenthaler2019overcoming,%
beck2020overcoming,%
hutzenthaler2020,%
hutzenthaler2020proof,%
HutzenthalerJentzenKruseNguyenWurstem2020,%
berner2020analysis}.

There are a number of criteria for classifying machine learning-based PDE solvers, among these mesh-free vs.\ fixed mesh methods, stochastic vs.\ deterministic methods or high-dimensional vs.\ low-dimensional methods.
While most of the investigated models can be considered mesh-free, we also mention some approaches that rely on an underlying and a priori known fixed mesh structure of the domain of the differential equations, cf.\ 
\cite{%
LeeKang1990,%
meade1994numerical,%
LagarisEtAl1998,%
LagarisEtAl2000,%
ramuhalli2005finite,%
malek2006numerical,%
chiaramonte2013solving,%
rudd2013solving}.

A method termed \emph{deep Galerkin method (DGM)} is proposed in \cite{SirignanoSpiliopoulos2018} and is applied to the solution of nonlinear second-order parabolic equations.
It is similar to the PINN approach discussed in~\cref{sec:PINNs} in that a neural network is used to approximate the PDE solution and the network is trained by minimizing a residual of the strong solution.
The methods is aimed at high-dimensional problems, however, and a Monte Carlo method rather than automatic differentiation is used to compute second derivatives. 
A similar approach for solving high-dimensional random PDEs by training a neural network on the strong or weak residual  is given in \cite{NabianMeidani2019}.
In \cite{berg2018unified} a deep neural network approximation to the solution of linear PDEs is constructed using the strong residual of the PDE as a loss function, similar to the PINN reviewed in~\cref{sec:PINNs}.

In~\cite{darbon2016algorithms}, an approach for solving a certain kind of high-dimensional first-order Hamilton-Jacobi equations is proposed based the Hopf formula \cite{hopf1965generalized} whose computational expense behaves polynomially in the spatial dimension.
In subsequent work, first-order Hamilton-Jacobi equations in high dimension are considered in~\cite{darbon2020overcoming,darbon2020some} based on classes of neural networks that exactly encode the viscosity solutions of these equations.%, similar to the technique presented in~\cref{sec:PINNs}.

Further approaches based on the multilevel decomposition of Picard approximations and on full-history recursive multilevel Picard approximations
\cite{%
eEtAl2019TR,%
eEtAl2017TR,%
giles2019generalised,%
eEtAl2019,%
hutzenthaler2020multilevel,%
beck2020overcoming} 
of type~\eqref{eq:SemilinPDE} have been successfully applied in high dimensions as well.
Other directions of research that deal with high-dimensional PDEs are branching diffusion processes \cite{henry2014numerical,henry2019branching}.

Another research area for the solution of PDEs is based on multi-scale \emph{deep neural networks (DNNs)}, cf.\ 
\cite{%
Cai2019MultiscaleDN,%
Liu2020MultiscaleDN,%
Wang_2020,%
Li_2020}.
In a recent \emph{Nature} publication multiscale DNNs are employed for diagnosing Alzheimer's disease \cite{lu2018multimodal}.
%\jb{<- Ich weiß nicht, ob das mit rein kann, habe ich gefunden und fand es spannend.}
Based on phase shift DNNs, \cite{Cai2020APS} considers the efficient solution of high-frequency wave equations.

In \cite{hure2021deep,bachouch2020deep} the authors introduce and compare a number of neural network-based algorithms applied to stochastic control problems, nonlinear PDEs and BSDEs, incl.\ the example discussed in \cref{sec:LinearQuadraticGaussianControl}.

While the approaches discussed so far employ neural networks to learn mappings between finite-dimensional Euclidean spaces, the methodologies proposed in
\cite{%
lu2019deeponet,%
bhattacharya2020model,%
nelsen2020random,%
li2020neural,%
li2020multipole,%
FourierNNs} 
aim to infer mappings between function spaces, known as \emph{neural operators}.
These mesh-free and infinite-dimensional operators require no prior knowledge of the underlying PDE but rely on a set of training data in the form of observations.

A general procedure based on data-driven machine learning to accelerate existing numerical methods for the solution of partial and ordinary differential equations is presented in \cite{mishra2018machine}.

A method to solve variational problems by means of scientific machine learning is proposed in~\cite{EYu2018}, termed the \emph{deep Ritz method} by the authors.
The method relies on a reformulation of variational problems as an energy minimization problem.
Boundary conditions are enforced weakly by the addition of a penalty term to the energy functional, for example
\begin{equation}
    \label{eq:RitzVar}
    \min_{u \in H} \int_\Omega \left( \frac{1}{2} \abs{\nabla u(x)}  - u(x)\right) \d x + \beta \, \int_{\partial \Omega} u(s)^2 \d s
\end{equation}
in the case of a Poisson problem with homogeneous boundary conditions, where $H$ is a set of admissible functions and $\beta$ is a penalty parameter used to enforce the boundary conditions.
The proposed methodology relies on three key ideas: the set of admissible functions $H$ is 
represented by a (deep) neural network; the integrals in the energy functional \eqref{eq:RitzVar} are approximated by Monte-Carlo sampling; and the neural network is trained through a stochatic gradient descent type algorithm on mini-batches.
An extension to this approach is given in \cite{liao2019deep} termed the \emph{deep Nitzsche method}.

Finally, we want to draw attention to the software package \texttt{NeuralPDE.jl} \cite{DifferentialEquations.jl-2017} written in the programming language Julia \cite{bezanson2017julia}.
It is available at \url{https://github.com/SciML/NeuralPDE.jl} and features the solution of PDEs by PINNs, forward-backward SDEs for parabolic PDEs as well as deep-learning based solvers for optimal stopping time problems and Kolmogorov backward equations.

\section{Conclusion} \label{sec:conclusion}

The methods reviewed in this paper illustrate the versatility of machine learning-based algorithms for the solution of PDEs and represent the currently most promising approaches.
While PINNs (\cref{sec:PINNs}) are, as of the writing of this survey, best suited for low-dimensional but complex nonlinear PDEs, the methods based on the Feynman-Kac theorem in \cref{sec:FeynmanKacSolver} and BSDEs in \cref{sec:DeepBSDESolver} promise to extend current simulation capabilities when employed for high-dimensional linear and semi-linear parabolic problems in non-variational form, for which classical approaches are infeasible due to the curse of dimensionality.
As deep learning continues to grow rapidly in terms of methodological, theoretical and algorithmic advances, we believe that the field of machine learning-based solution methods of PDEs promises to remain an exciting research field in the  years ahead.

\bibliography{survey_HJBs.bib}
\end{document}